\documentclass[3p]{elsarticle}
\usepackage[utf8]{inputenc}
\usepackage{xcolor}
\usepackage{amsmath, amsthm, amssymb, amsfonts, mathtools}
\usepackage[colorlinks=true]{hyperref}
\usepackage{graphicx}
\usepackage{caption}
\usepackage{subcaption}
\usepackage{tikz}

\theoremstyle{plain}
\newtheorem{theorem}{Theorem}[section]
\newtheorem{proposition}{Proposition}[section]

\theoremstyle{remark}
\newtheorem{remark}{Remark}
\theoremstyle{definition}
\newtheorem{definition}{Definition}[section]

\begin{document}
\begin{frontmatter}
\title{Strategic control for a Boltzmann-like decision-making model}
\author[1]{Luis Guillermo Venegas-Pineda}
\author[1]{Hildeberto Jard\'on-Kojakhmetov\corref{cor1}}
\ead{h.jardon.kojakhmetov@rug.nl}
\author[2,3]{Maximilian Engel}
\author[4]{Jobst Heitzig}
\author[2]{Muhittin Cenk Eser}
\author[5]{Ming Cao}
\cortext[cor1]{Corresponding author}
\affiliation[1]{organization={Bernoulli Institute for Mathematics, Computer Science and Artificial Intelligence, University of Groningen, Nijenborgh 9, 9700 AK, Groningen, The Netherlands.}}
\affiliation[2]{organization={Institute of Mathematics, Freie Universit{\"a}t Berlin, Arnimallee 6, 14195 Berlin, Germany.}}
\affiliation[3]{organization={Korteweg-de Vries Institute for Mathematics, University of Amsterdam, Science Park 105-107, 1098 XG Amsterdam, The Netherlands.}}
\affiliation[4]{organization={Complexity Science Department, Potsdam Institute for Climate Impact Research, PO Box 60 12 03, 14412 Potsdam, Germany.}}
\affiliation[5]{organization={Engineering and Technology Institute Groningen, University of Groningen, Nijenborgh 9, 9700 AE, Groningen, The Netherlands.}}

\begin{abstract}
 We study a prototypical non--polynomial decision--making model for which agents in a population potentially alternate between two consumption strategies, one related to the exploitation of an unlimited but considerably expensive resource and the other a comparably cheaper but restricted and slowly renewable source. In particular, we study a model following a Boltzmann--like exploration policy, enhancing the accuracy at which the exchange rates are captured with respect to classical polynomial approaches by considering sigmoidal functions to represent the cost--profit relation in both exploit strategies. Additionally, given the intrinsic timescale separation between the decision--making process and recovery rates of the renewable resource, we use geometric singular perturbation theory to analyze the model. We further use numerical analysis to determine parameter ranges for which the model has a distinct number of fold points of its critical manifold. These points, being related to critical states of the system, are relevant to the fast transitions between strategies. Hence, we design controllers to regulate such rapid transitions by taking advantage of the system's criticality.
\end{abstract}
\end{frontmatter}

\section{Introduction}\label{Sec:Introduction}
A recurrent problem across the sciences is decision--making, which can be summarized as a rational agent
choosing among alternatives so as to maximise some expected profit
\cite{Ref:Kelly2003, Ref:Peterson2009}, on time scales ranging from operational to
strategic and political \cite{Ref:Marugan2015, Ref:GarciaMarquez2022}, with
applications from biology \cite{Ref:Helikar2008, Ref:Balazsi2011} to economics
\cite{Ref:Lynn2015} and psychology \cite{Ref:Mishra2014}. A standard way to
capture the resulting bounded rational behaviour is through exponential weighting
\cite{Ref:Arora2012}, which in reinforcement learning \cite{Ref:Sutton2018, Ref:Wang2020}
appear as softmax, Gibbs, or \emph{Boltzmann} exploration policies, in which
the probability of selecting a strategy is proportional to an exponential function
of its expected reward and which balance exploitation against exploration
\cite{Ref:Bianchi2017}.

    For our research, we study a resource's stock dynamics as being consumed by two distinguished groups following different exploitation strategies, one consuming an unlimited but highly costly common resource, such as wind energy, and the second employing a comparably cheaper but restricted and slowly renewable resource, for instance biomass. Let $y \geq 0$ represent the limited resource stock and $x \in[0,1]$ the share of agents exploiting it, while $1-x$ is the portion of agents consuming the unlimited resource instead. 
    
    We model the joint evolution of $x$ and $y$ as the fast--slow system
    \begin{subequations}
    \renewcommand{\theequation}{\theparentequation.\arabic{equation}}
        \begin{align}
            \dot{x} &= \gamma_{1}(1-x)\left( \eta_{1} + \frac{1-\eta_{1}}{1+e^{-\beta_{1}(\alpha_{1} + \delta(x,y))}} \right) - \gamma_{2} x \left( \eta_{2} + \frac{1-\eta_{2}}{1+e^{-\beta_{2}(\alpha_{2} - \delta(x,y))}} \right), \label{Eq:FastComponentOriginalSystem}\\
            \dot{y} &= \varepsilon y (1-rx), \label{Eq:SlowComponentOriginalSystem}
        \end{align}
        \label{Eq:OriginalSystem}
    \end{subequations}\par
    \noindent where $0<\varepsilon\ll1$ indicates the timescale separation, involving the \emph{slow} recovery speed of the limited resource stock $y$ and the comparatively fast change of $x$ due to the agents' adaptation of their exploitation strategies. 
    
    \begin{remark}
        Much of our analysis relies on (standard techniques of) geometric singular perturbation theory (GSPT). For completeness, we include the necessary background on GSPT in \ref{Sec:Preliminaries}, complemented by relevant references. 
    \end{remark}
    
    The \emph{natural} and \emph{metabolic} component of \eqref{Eq:OriginalSystem} is a simple equation for $y$ governed by the growth rate $\varepsilon$ and each agents' relative harvesting rate $r$, resulting in an effective total harvesting rate of $\varepsilon r x$. On the other hand, the \emph{economic} element is a model of bounded rational behaviour governed by a set of parameters as follows. The terms $\gamma_{1,2}>0$ represent the rate at which agents from one strategy consider switching to the opposite strategy. If an agent considers switching, they either change strategy independently of the possible profits, which happens with probability $\eta_{1,2}\in[0,1]$ and can be called \emph{unconditional exploration}, or they base their decision whether to switch strategies on the profit difference 
    \begin{equation}
        \delta(x,y) = y + \frac{c}{d-x}-b,
        \label{Eq:ProfitDifference}
    \end{equation}
    which happens with probability $1-\eta_{1,2}$, and where $b>1$ represents the benefits of harvesting one unit of the unlimited resource, while $c>0$ and $d>1$ are parameters governing the costs of harvesting one unit of the unlimited resource. Note that the function $\delta(x,y)$ is chosen such that these costs decrease from $c/(d-1)$ to $c/d$ as the share of agents exploiting it, i.e.~$1-x$, grows. Particularly, agents following the latter strategy employ a Boltzmann (or softmax) policy to decide their future strategy \cite{Ref:Achbany2008, Ref:Bianchi2017}. Therefore, the logits, i.e., the logarithm of the odds of the probability of a certain event occurring, of the relative probabilities for both strategies are proportional to the corresponding profits, multiplied by an \emph{inverse temperature} parameter $\beta_{1,2}>0$, which is a common way to model bounded rationality \cite{Ref:Herbert1997, Ref:Lorkowski2017}. In addition to the payoff comparison, agents favour switching to a certain extent governed by an \emph{offset} parameter $\alpha_{1,2}>0$, which can be interpreted as a \emph{conditional exploration}. Altogether, considering soft optimization and conditional exploration lead to the logistic terms $1 / (1 + e^{-\beta(\alpha\pm\delta)})$ of \eqref{Eq:OriginalSystem}.  
    Similar approaches are found in economics, concerning the aggregate saving rate in a large population \cite{Ref:Asano2021}, as well as in ecology, where the cost--profit difference between distinct energy resources is expressed \cite{Ref:Zeppini2020} or the rate of succession of grassland to forest as described in
    \cite{Ref:Innes2013}. 

Beyond reinforcement learning, exponential and logit choice rules have a
long history in economics and learning theory
\cite{mcfadden1972conditional,sutton1998reinforcement}. In parallel, fast--slow dynamical systems exhibiting fold singularities and canard dynamics have been extensively studied from a mathematical perspective \cite{Ref:Krupa2001,Ref:Kuehn2015} and in a variety of applied contexts. 

    A representative recent case study is the tutorial treatment of the
FitzHugh--Nagumo system~\cite{goncalves},
which combines a codimension-one bifurcation analysis with a numerical
determination of the canard locus in parameter space. Our fold-point diagrams
(Figure~\ref{fig:fold-diagrams}) and singular-Hopf criticality maps
(Figure~\ref{fig:criticalitymaps}) play an analogous role for the present
non-polynomial Boltzmann-type model, identifying where folds, singular Hopf
bifurcations and canard explosions occur as the cost parameters vary.

\begin{figure*}[htbp]
    \centering
        \begin{subfigure}[htbp]{0.5\textwidth}
        \includegraphics[]{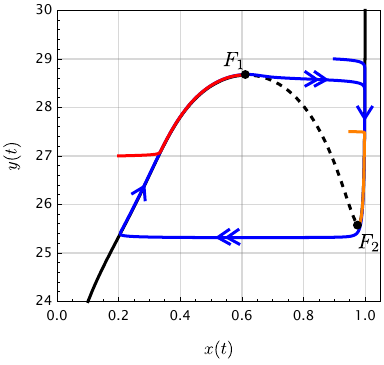}
        \end{subfigure}\hfill
        \begin{subfigure}[htbp]{0.5\textwidth}
        \includegraphics[]{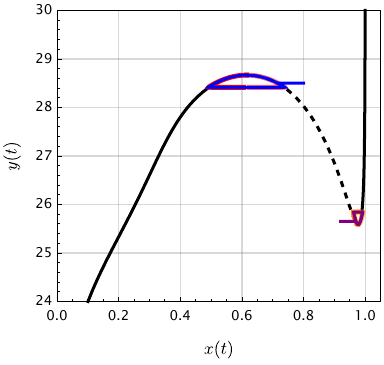}
        \end{subfigure}
        \caption{Open--loop (left) and controlled (right) responses of the decision--making model \eqref{Eq:OriginalSystem}, where the stable and unstable regions of the critical manifold, see \ref{Sec:FastSlowSystems}, are depicted in solid and dashed black, respectively. For the open--loop (left), a relaxation oscillation (blue) and two trajectories (red, orange) stabilizing at the closest fold point ($F_1, F_2$ respectively) are shown. For the closed-loop (right), we present the stabilization of two trajectories (blue, purple) of \eqref{Eq:OriginalSystem} along the desired reference orbits (thick red, orange) through the implementation of the controlled schemes designed in this work, see Section \ref{Sec:Results}. Parameters: $\alpha = 2.0$, $\beta = 0.75$, $\gamma = 0.5$, $c = 2.5$, $d = 1.18$, $b = 30.0$, $\varepsilon = 0.001$. For the open--loop scenario: relaxation oscillation $\left( x(0), y(0) \right) = (0.9, 29.0)$ with $r = 1.62$, trajectory stabilizing at $F_1$ $\left( x(0), y(0) \right) = (0.2, 27.0)$ with $r = 1.65$, and trajectory stabilizing at $F_2$ $\left( x(0), y(0) \right) = (0.95, 27.5)$ with $r = 1.08$. Double/single arrows mark fast/slow segments.}
        \label{Fig:CompetitiveSystemOpenLoop}
    \end{figure*}

    For generic initial conditions, the flow generated by \eqref{Eq:OriginalSystem} quickly converges to attracting equilibria of the fast dynamics \eqref{Eq:FastComponentOriginalSystem}. Once a trajectory is sufficiently close to such equilibria, it evolves according to the slow dynamics \eqref{Eq:SlowComponentOriginalSystem} until the equilibria of the fast dynamics undergo a bifurcation, repeating the aforementioned process to either produce sustained oscillations or stabilize in an equilibrium point of \eqref{Eq:OriginalSystem}. For more details on fast--slow systems, see \ref{Sec:FastSlowSystems}. In particular, our numerical results show that under parameter variations, the critical manifold of system \eqref{Eq:OriginalSystem} presents either zero, two or four fold points, see Section \ref{Sec:Assumptions}. The previous fact provides our model with the capability to produce interesting dynamics, such as those shown in Figure \ref{Fig:CompetitiveSystemOpenLoop} for a case with two fold points. In the left panel, open--loop responses of system \eqref{Eq:OriginalSystem} are depicted, including the stabilization of trajectories (red, orange) to equilibrium points near the fold points $F_{1,2}$, and a relaxation oscillation (blue), a response formed by alternating fast and slow segments \cite{Ref:Grasman2011}. On the other hand, the right panel exhibits actuated responses (blue, purple) of \eqref{Eq:OriginalSystem} for the stabilization of \emph{canard orbits} (red, orange), particular solutions known for their lack of robustness towards perturbations as well as their capability to travel for considerable amounts of time along repelling regions of the critical manifold, geometric object associated to \eqref{Eq:OriginalSystem} which describes the boundary between the fast \eqref{Eq:FastComponentOriginalSystem} and slow \eqref{Eq:SlowComponentOriginalSystem} dynamics when $\varepsilon=0$ \cite{Ref:Szmolyan2001, Ref:Kuehn2015, Ref:JardonKojakhmetov2022}.

We now state the motivation for our work and highlight its novelty.

    \emph{Motivation:} In the uncontrolled system, the canard cycles
associated with the fold points $F_{1,2}$ are structurally unstable, i.e., they exist
only within an exponentially small window of the harvesting rate $r$ and are
destroyed by arbitrarily small perturbations, so the criticality that drives the
abrupt transitions between strategies cannot be exploited reliably. This motivates
a control approach that renders a chosen canard cycle a robust attractor, enabling
a controlling entity to exploit the system's criticality, and hence to produce or
suppress drastic consumption changes, with minor effort. From an applied
perspective, stabilising an orbit near $F_1$ amounts to preserving a high stock of
the limited resource $y$, for instance a public good such as a forest providing
ecological habitats or recreational value, by limiting the share of agents
consuming it. In contrast, stabilising near $F_2$ amounts to consuming $y$ in much
larger but still sustainable amounts for purely economic reasons.

    \emph{Novelty:} Existing studies of such decision dynamics
typically rely on polynomial or locally truncated models, and control objectives
are seldom posed near non--generic folds. Our contribution is threefold. First, we
analyse a genuinely non--polynomial, Boltzmann--type decision model within a
fast--slow framework, and we locate its folds together with the associated singular
Hopf bifurcations and canard explosions (Section~\ref{Sec:Assumptions}). Second,
and central to this work, we design compatible fast--slow controllers, that is,
schemes that stabilise a target canard cycle without dramatically altering the
original dynamical structure (see \ref{Sec:PlanarFoldedCanards}), around a possibly
non--generic fold point. This extends the generic theory of
\cite{Ref:JardonKojakhmetov2022} to higher contact orders $k>1$. The relevance of non-generic folds, where higher-order terms better capture the local geometry of the critical manifold, is detailed in Section~\ref{Sec:Results}. Third, the
proposed controllers are robust to bounded parametric uncertainty and act on the
agents' decision rules rather than on the resource dynamics directly, which makes
them less intrusive, and more realistic, from an applied standpoint.

    The rest of this work is organized as follows. In Section \ref{Sec:Assumptions} we present and justify some assumptions that simplify our analysis by reducing the effective number of parameters to be considered. Next, our main results are detailed in Section \ref{Sec:Results}. First, we derive compatible fast--slow controllers for a generalized quadratic system, see Theorem \ref{Thm:GeneralizedQuadraticControl}. Thereafter, we give arguments for applying our results in the decision--making system \eqref{Eq:OriginalSystem} and demonstrate its effectiveness to control canard cycles, even in the presence of modelling uncertainties. Finally, we discuss the reach of our controllers and potential future study directions in Section \ref{Sec:Discussion}. In \ref{Sec:Preliminaries}, we present the mathematical background necessary for our studies.

\section{Preliminary analysis and motivation}\label{Sec:Assumptions}
System \eqref{Eq:OriginalSystem} has thirteen parameters $(\gamma_1,\gamma_2,\eta_1,\eta_2,\alpha_1,\alpha_2,\beta_1,\beta_2,b,c,d,r,\varepsilon)$, which adds to the complexity of its analysis. To mitigate this complication, we describe and justify the main assumptions we make about \eqref{Eq:OriginalSystem} to reduce the number of parameters. In addition, we provide a motivation for our main results, namely, the control of `small' canard cycles.
The following assumptions are made for \eqref{Eq:OriginalSystem}. 

\begin{description}
    \item[A1:]  The parameters $\eta_i$
 represent the probability that agents switch strategy \emph{independently of payoffs}. In the context considered here, namely switching between exploiting a slowly renewable but cheap resource and an unlimited but costly resource, such payoff-independent switches are expected to be rare, since strategic changes are mainly driven by relative profitability. It is therefore natural to assume $\eta_i$ is small. Hence, for purposes of analysis, we assume that $\eta_1=\eta_2=0$.
 \item[A2:] The offset parameters $\alpha_i$ represent conditional exploration biases, independent of the actual payoff difference. In the absence of compelling evidence for an asymmetric bias, and to reduce parameter dimensionality, we set $\alpha_1=\alpha_2=\alpha>0$. 
 \item[A3:] The inverse temperature parameters $\beta_i$ model how strictly agents react to payoff differences. Since both groups are drawn from the same population of decision-makers and face similar informational constraints, it is reasonable to assume a common value $\beta_1=\beta_2=\beta>0$.
\end{description}

We note that the qualitative features of our results remain valid under small perturbations of A1–A3. In fact, all our simulations are performed without the above considerations to verify such a claim.

Next, let $\gamma\coloneqq\frac{\gamma_1}{\gamma_2}$. After a time rescaling $t\mapsto \frac{t}{\gamma_2}$ and recycling the notation $\varepsilon$ for $\frac{\varepsilon}{\gamma_2}$ we obtain the simplified model
\begin{equation}
    \begin{split}
        \dot{x} &=   \frac{\gamma(1-x)}{1+e^{-\beta(\alpha + \delta(x,y))}}  -    \frac{x}{1+e^{-\beta(\alpha - \delta(x,y))}}, \\
            \dot{y} &= \varepsilon y (1-rx).
    \end{split}
\end{equation}
Finally, we perform the rescaling  $(\tilde\alpha,\tilde y, \tilde c, \tilde b)=(\beta\alpha,\beta y,\beta c,\beta b)$ to eliminate the parameter $\beta$. For notational convenience we recycle the notation $\alpha,y,c,b$ for the rescaled terms obtaining
\begin{equation}\label{eq:main2}
    \begin{split}
        \dot{x} &=   \frac{\gamma(1-x)}{1+e^{-(\alpha + \delta(x,y))}}  -    \frac{x}{1+e^{-(\alpha - \delta(x,y))}}, \\
            \dot{y} &= \varepsilon y (1-rx).
    \end{split}
\end{equation}
With this, we have reduced the number of parameters to seven $(\gamma,\alpha,b,c,d,r,\varepsilon)$.

 The critical manifold of \eqref{eq:main2} (see the necessary background on fast-slow systems in \ref{Sec:Preliminaries}) is defined as
\begin{equation}\label{eq:C0}
    \mathcal{C}_0=\left\{ (x,y)\in\mathbb R^2\,|\frac{\gamma(1-x)}{1+e^{-(\alpha + \delta(x,y))}}  -    \frac{x}{1+e^{-(\alpha - \delta(x,y))}}=0\,\right\}.
\end{equation}
This already highlights the difficulty of analyzing \eqref{eq:main2}, as the critical manifold does not have a clear explicit expression, with the exponential function playing the most important role. Hence, we combine analytical and numerical methods to describe $\mathcal{C}_0$. The following propositions provides certain sufficient conditions describing the properties of the critical manifold, and is complemented numerically afterwards.
\begin{proposition}\label{prop:NH_simple}
Assume $\alpha>0$, $\gamma>0$, $c>0$, and $d>1$. If
\[
c\le 4d(d-1),
\]
then the critical manifold $\mathcal C_0$ is normally hyperbolic on $x\in[0,1]$.

\end{proposition}

\begin{proof}
Let
\[
f(x,y)=\gamma(1-x)\sigma(\alpha+\delta(x,y))-x\,\sigma(\alpha-\delta(x,y)),
\qquad
\delta(x,y)=y+\frac{c}{d-x}-b,
\]
where $\sigma(z)=(1+e^{-z})^{-1}$. The critical manifold is
\[
\mathcal C_0=\{(x,y)\in\mathbb R^2:f(x,y)=0\}.
\]
Since $x$ is the fast variable, loss of normal hyperbolicity can occur only at points of $\mathcal C_0$ where $\partial_x f=0$.

Let
\[
A=\sigma(\alpha+\delta(x,y)),
\qquad
B=\sigma(\alpha-\delta(x,y)).
\]
Then
\[
f(x,y)=\gamma(1-x)A-xB.
\]
On $\mathcal C_0$ one has
\[
\gamma(1-x)A=xB,
\]
and the condition $\partial_x f=0$ reduces to
\[
\frac{c\,x(1-x)}{(d-x)^2}(2-A-B)=1.
\]
Since $\alpha>0$,
\[
A+B
=
1+\frac{\sinh\alpha}{\cosh\alpha+\cosh\delta}
>1,
\]
hence
\[
2-A-B<1.
\]
Therefore, any nonhyperbolic point must satisfy
\begin{equation}\label{eq:cond_nh}
    1<\frac{c\,x(1-x)}{(d-x)^2}.
\end{equation}

Now define
\[
m(x)=\frac{x(1-x)}{(d-x)^2},
\qquad x\in[0,1].
\]
A direct computation shows that $m$ attains its maximum at
\[
x_*=\frac{d}{2d-1},
\]
and
\[
\max_{x\in[0,1]}m(x)=\frac{1}{4d(d-1)}.
\]
Hence loss of normal hyperbolicity requires
\[
1<\frac{c}{4d(d-1)}.
\]
The claim follows.
\end{proof}

Intuitively, normal hyperbolicity asks that the fast flow be uniformly attracting
or repelling transverse to $\mathcal{C}_0$. Equivalently, that the fast
linearisation at $\mathcal{C}_0$ have no eigenvalue on the imaginary axis. For the scalar fast variable
here this is simply $\partial_x f\neq0$ ($\partial_x f<0$ attracting,
$\partial_x f>0$ repelling), and it fails precisely at the folds, where
$\partial_x f=0$. The formal definition is given in
\ref{Sec:FastSlowSystems} (Definition~\ref{def:nh}). Contrasting Proposition \ref{prop:NH_simple}, the next Proposition qualitatively characterizes the non-hyperbolic scenario.

\begin{proposition}\label{prop:fold_pairs}
Assume $\alpha>0$, $\gamma>0$, $c>0$, and $d>1$. If the critical manifold $\mathcal C_0$ is not normally hyperbolic, then, generically, its nonhyperbolic points are folds and they occur in an even number. In particular, generically, $\mathcal C_0$ has at least two fold points.
\end{proposition}

\begin{proof}
Define
\[
h(x)=\frac{c}{d-x},
\qquad
u=y+h(x)-b,
\qquad
A(u)=\sigma(\alpha+u),
\qquad
B(u)=\sigma(\alpha-u).
\]
The critical manifold equation can be written as
\[
\gamma(1-x)A(u)-xB(u)=0,
\]
that is,
\[
x=R(u):=\frac{\gamma A(u)}{\gamma A(u)+B(u)}.
\]
Since
\[
A'(u)=A(u)(1-A(u)),
\qquad
B'(u)=-B(u)(1-B(u)),
\]
one finds
\[
R'(u)=\frac{\gamma A(u)B(u)(2-A(u)-B(u))}{(\gamma A(u)+B(u))^2}>0.
\]
Thus $R$ is strictly increasing, with
\[
\lim_{u\to-\infty}R(u)=0,
\qquad
\lim_{u\to\infty}R(u)=1.
\]
Hence $\mathcal C_0$ admits the parametrization
\[
\phi:u\mapsto (R(u),u-h(R(u))+b).
\]

Fold points correspond to points where the tangent of this parametrized curve is vertical, namely to solutions of
\[
1-h'(R(u))R'(u)=0.
\]
Equivalently, writing
\[
\Phi(u):=h'(R(u))R'(u),
\]
one has
\begin{equation}\label{eq:Phi}
    \Phi(u)=\frac{c\gamma A(u)B(u)(2-A(u)-B(u))}{(dB(u)+(d-1)\gamma A(u))^2}.
\end{equation}

Since $A(u),B(u)\in(0,1)$ and
\[
\lim_{u\to\pm\infty}A(u)B(u)=0,
\]
it follows that
\[
\lim_{u\to\pm\infty}\Phi(u)=0.
\]
Therefore, the equation $\Phi(u)=1$ has, generically, an even number of solutions. If $\mathcal C_0$ is not normally hyperbolic, then $\Phi(u)=1$ has at least one solution, and hence, generically, at least two. For generic parameter values these nonhyperbolic points are of fold type.
\end{proof}

\begin{figure}[htbp]
    \centering
    \begin{tikzpicture}
        \node at (0,0) {
        \begin{tikzpicture}
            \node at (0,0){
            \includegraphics[]{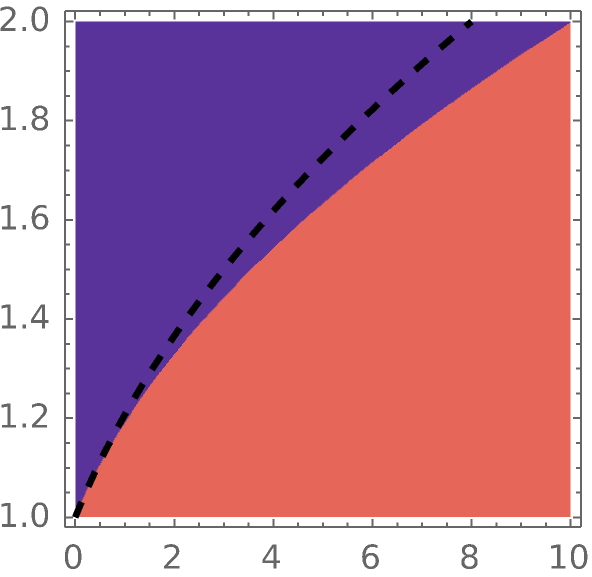}
            };
            \node at (0,-2.75){$c$};
            \node at (-3,0){$d$};
            \node at (0,2.75){$\alpha=\frac{1}{2}\ln3$, $\gamma=1$};
        \end{tikzpicture}
        };
        \node at (6,0) {
        \begin{tikzpicture}
            \node at (0,0){
            \includegraphics[]{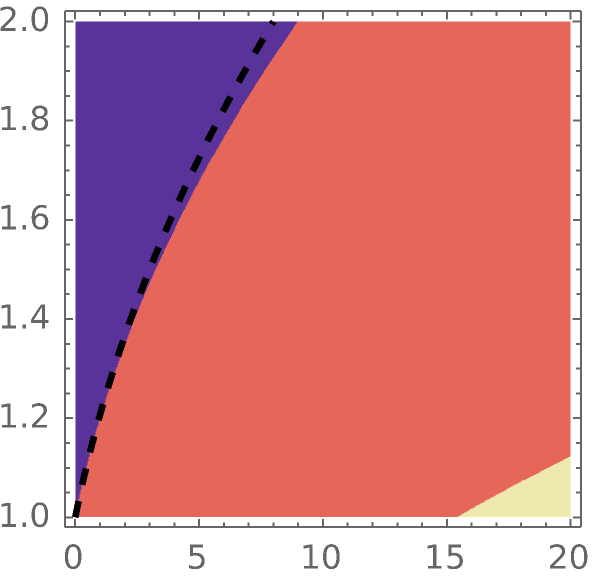}
            };
            \node at (0,-2.75){$c$};
            \node at (-3,0){$d$};
            \node at (0,2.75){$\alpha=2\ln3$, $\gamma=\frac{1}{4}$};
        \end{tikzpicture}
        };
        \node at (0,-6.5) {
        \begin{tikzpicture}
            \node at (0,0){
            \includegraphics[]{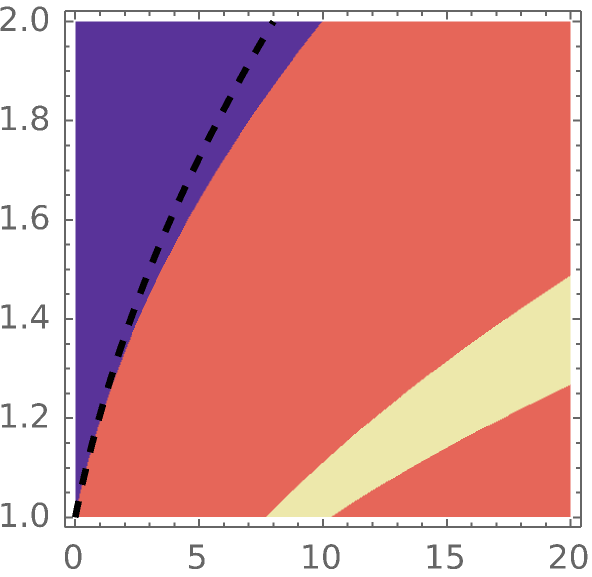}
            };
            \node at (0,-2.75){$c$};
            \node at (-3,0){$d$};
            \node at (0,2.75){$\alpha=2\ln3$, $\gamma=\frac{1}{2}$};
        \end{tikzpicture}
        };
        \node at (6,-6.5) {
        \begin{tikzpicture}
            \node at (0,0){
            \includegraphics[]{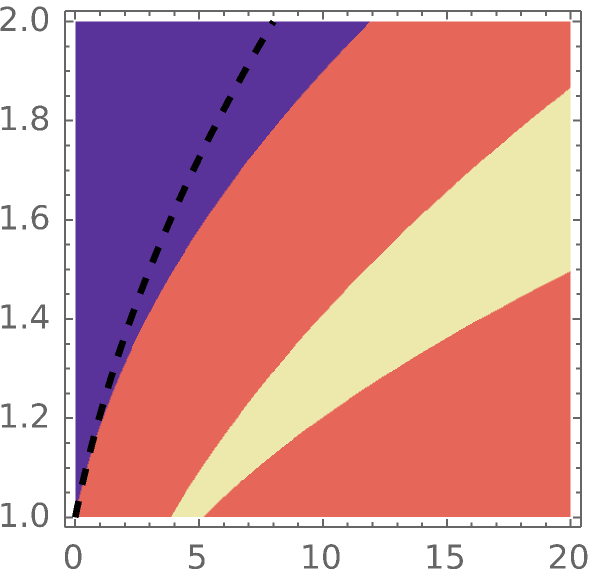}
            };
            \node at (0,-2.75){$c$};
            \node at (-3,0){$d$};
            \node at (0,2.75){$\alpha=2\ln3$, $\gamma=1$};
        \end{tikzpicture}
        };
        \node at (0,-13) {
        \begin{tikzpicture}
            \node at (0,0){
            \includegraphics[]{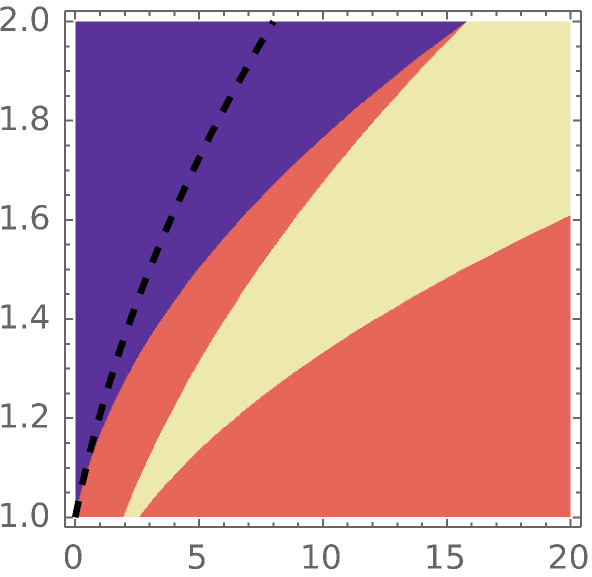}
            };
            \node at (0,-2.75){$c$};
            \node at (-3,0){$d$};
            \node at (0,2.75){$\alpha=2\ln3$, $\gamma=2$};
        \end{tikzpicture}
        };
        \node at (6,-13) {
        \begin{tikzpicture}
            \node at (0,0){
            \includegraphics[]{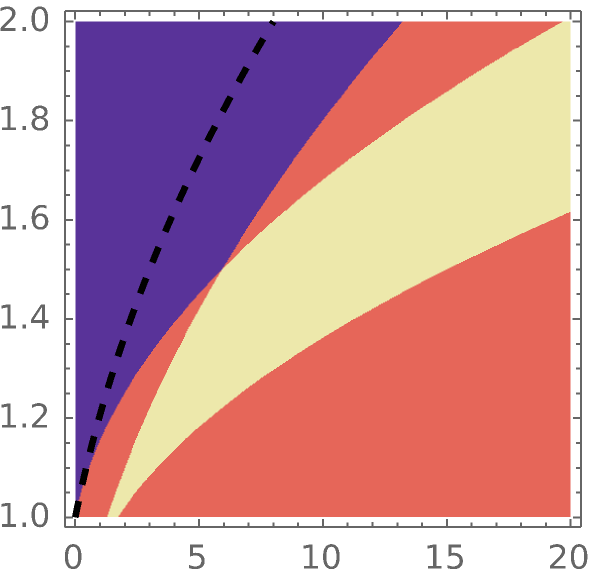}
            };
            \node at (0,-2.75){$c$};
            \node at (-3,0){$d$};
            \node at (0,2.75){$\alpha=2\ln3$, $\gamma=3$};
        \end{tikzpicture}
        };

    \end{tikzpicture}
 
    \caption{Fold-point diagrams representing the number of fold points the critical manifold of \eqref{eq:main2} has as a function of the parameters $(c,d)$ and with the parameters $\alpha$ and $\gamma$ fixed as indicated. Colours in each diagram correspond to the number of fold points: purple, red, and yellow represent $0$, $2$, and $4$ fold points, respectively. See also Figure \ref{fig:criticalmanifolds}. The dashed curve represents the sufficient condition of Proposition \ref{prop:NH_simple}. Due to the fact that \eqref{eq:main2} is a fast-slow system, these fold points also indicate singular Hopf bifurcations; their criticality is further discussed in Proposition \ref{prop:hopf}, see also Figure \ref{fig:criticalitymaps}.}
    \label{fig:fold-diagrams}
\end{figure}

\begin{figure}
    \centering
    \begin{tikzpicture}
        \node at (0,0){
        \begin{tikzpicture}
            \node at (0,0){
            \includegraphics[scale=1.2]{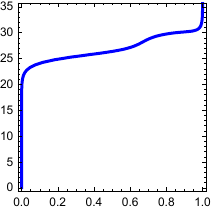}
            };
            \node at (0,-2.5){$x$};
            \node at (-2.5,0){$y$};
            \node at (0,2.5){$c=1$};
            \end{tikzpicture}
        };
        \node at (6,0){
        \begin{tikzpicture}
            \node at (0,0){
            \includegraphics[scale=1.2]{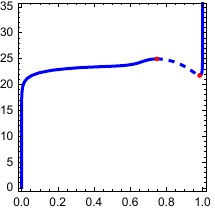}
            };
            \node at (0,-2.5){$x$};
            \node at (-2.5,0){$y$};
            \node at (0,2.5){$c=3$};
            \end{tikzpicture}
        };
        \node at (0,-6){
        \begin{tikzpicture}
            \node at (0,0){
            \includegraphics[scale=1.2]{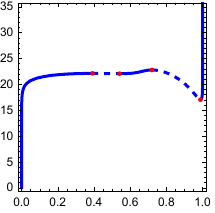}
            };
            \node at (0,-2.5){$x$};
            \node at (-2.5,0){$y$};
            \node at (0,2.5){$c=4$};
            \end{tikzpicture}
        };
        \node at (6,-6){
        \begin{tikzpicture}
            \node at (0,0){
            \includegraphics[scale=1.2]{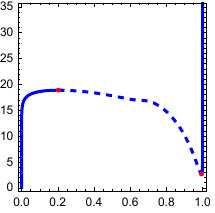}
            };
            \node at (0,-2.5){$x$};
            \node at (-2.5,0){$y$};
            \node at (0,2.5){$c=7$};
            \end{tikzpicture}
        };
    \end{tikzpicture}
    \caption{Critical manifold $\mathcal C_0$ of \eqref{eq:main2} for four values of the cost
parameter $c$ (others fixed), illustrating the transition through $0$, $2$ and $4$ fold
points. Attracting and repelling branches are drawn solid and dashed; red dots mark the
fold points. Parameters $\gamma=2,\ d=1.2,\ \alpha=2\ln 3,\ b=30$;
compare with the lower-left panel of Figure~\ref{fig:fold-diagrams}.
As $c$ varies, folds are created or annihilated in pairs at fold collisions; near
such a collision the colliding folds are closely spaced and individually flat, so locally
the manifold is better captured by a higher-order than by the generic quadratic
approximation. The four-fold panel ($c=4$) illustrates this regime, further discussed in
Section~\ref{Sec:Applications}.}
    \label{fig:criticalmanifolds}
\end{figure}

\begin{remark}\label{rem:parameter_planes}
To further explore the properties of the critical manifold beyond the previous results, we present ``fold-point diagrams'' in parameter space. The parameter planes are chosen according to the roles of the parameters in the fold condition. The primary plane is $(c,d)$, since the explicit bound
\[
c\le 4d(d-1)
\]
depends only on these parameters, and they govern the factor (see \eqref{eq:cond_nh})
\[
\frac{c\,x(1-x)}{(d-x)^2}.
\]
The parameters $\alpha$ and $\gamma$ modulate this primary picture in different ways. To see this, let $N$ denote the numerator of $\Phi$ \eqref{eq:Phi}. The parameter $\alpha$ affects the sigmoidal response and, in particular, the change of type at $u=0$ is detected by
\[
N''(0)=\frac{2c\gamma e^{2\alpha}(3-e^\alpha)}{(1+e^\alpha)^5},
\]
which vanishes precisely at $\alpha=\ln 3$. This motivates distinguishing between the cases $\alpha<\ln 3$ and $\alpha>\ln 3$ in the $(c,d)$ fold diagrams shown in Figure \ref{fig:fold-diagrams}; a sample of critical manifolds is further shown in Figure \ref{fig:criticalmanifolds}.

The parameter $\gamma$ enters through the balance relation on the critical manifold. Although it does not appear in the simple exclusion bound above, it changes the shape of the fold locus in the $(c,d)$-plane in a nontrivial way, especially for $\alpha>\ln 3$. For this reason, the $(c,d)$ diagrams are presented for representative fixed values of $\gamma$.

By contrast, the parameter $b$ only shifts the critical manifold in the $y$-direction and is therefore not used in the diagrams.
\end{remark}

We close this section by locating the bifurcations of the full system
$0<\varepsilon\ll 1$ that organize the transition between equilibria,
canard cycles, and relaxation oscillations. Although the critical
manifold $\mathcal{C}_0$ does not admit a closed-form expression, the
structure of the slow equation in \eqref{eq:main2} allows us to
determine the location of the Hopf bifurcations \emph{exactly}.
Throughout, let
\begin{equation*}
f(x,y)=\gamma(1-x)\,\sigma\big(\alpha+\delta(x,y)\big)
-x\,\sigma\big(\alpha-\delta(x,y)\big)
\end{equation*}
denote the right-hand side of the fast equation in \eqref{eq:main2},
where $\sigma(z)=\left(1+e^{-z}\right)^{-1}$.
\begin{proposition}
\label{prop:hopf}
Let $r>1$ and let $y^{*}(r)$ denote the unique solution of
$f(1/r,y)=0$, and assume $y^{*}(r)>0$. Then the equilibrium
$p(r)=\left(1/r,\,y^{*}(r)\right)$ of \eqref{eq:main2} satisfies, for
every $\varepsilon>0$,
\begin{equation*}
\operatorname{tr} J(p(r))=\partial_x f(p(r)),
\qquad
\det J(p(r))=\varepsilon\, r\, y^{*}(r)\,\partial_y f(p(r))>0,
\end{equation*}
where $J$ denotes the Jacobian of \eqref{eq:main2}. In particular:
\begin{enumerate}
\item The equilibrium point $p(r)$ is never a saddle, and its stability
coincides with the normal stability of the branch of $\mathcal{C}_0$ on
which it lies: $p(r)$ is attracting if $\partial_x f(p(r))<0$ and
repelling if $\partial_x f(p(r))>0$.
\item If $F=(x_F,y_F)$ is a fold point of $\mathcal{C}_0$, then
\eqref{eq:main2} undergoes a Hopf bifurcation at
\begin{equation*}
r_H=\frac{1}{x_F},
\end{equation*}
exactly and independently of $\varepsilon$, with eigenvalues
$\pm i\omega$, where
$\omega=\sqrt{\varepsilon\, r_H\, y_F\, \partial_y f(F)}
=\mathcal{O}(\sqrt{\varepsilon})$, so that the bifurcation is singular
\cite{Ref:Krupa2001}.
\item The eigenvalue crossing is transversal if and only if
$\partial_{xx}f(F)\neq 0$, i.e., if and only if the fold is generic
($k=1$); in that case
$\frac{d}{dr}\operatorname{tr}J(p(r))\big\vert_{r=r_H}
=-\partial_{xx}f(F)/r_H^{2}$.
\end{enumerate}
\end{proposition}
\begin{proof}
Recall that $\sigma'(z)=\sigma(z)\left(1-\sigma(z)\right)>0$ for all
$z\in\mathbb{R}$. Hence, one has
\begin{equation*}
\partial_y f(x,y)=\gamma(1-x)\,\sigma'(\alpha+\delta(x,y))
+x\,\sigma'(\alpha-\delta(x,y))>0
\qquad \text{for all } x\in(0,1).
\end{equation*}
Moreover, $f(x,y)\to\gamma(1-x)>0$ as $y\to\infty$ and
$f(x,y)\to -x<0$ as $y\to-\infty$, so $y^{*}(r)$ exists and is unique
for every $x^{*}=1/r\in(0,1)$; in particular, $\mathcal{C}_0$ intersects
each vertical line $\{x=\mathrm{const}\}$, $x\in(0,1)$, exactly once.
The Jacobian of \eqref{eq:main2} reads
\begin{equation*}
J(x,y)=\begin{pmatrix}
\partial_x f & \partial_y f\\[2pt]
-\varepsilon r y & \varepsilon(1-rx)
\end{pmatrix},
\end{equation*}
and at $p(r)$ the entry $\varepsilon(1-rx)$ vanishes, which yields the
stated expressions for the trace and the determinant. Item~1 follows
immediately, since $\det J>0$ excludes saddles and the trace equals
$\partial_x f$, whose sign determines the normal stability of
$\mathcal{C}_0$. For item~2, note that $p(r)\in\mathcal{C}_0$ for every
$r$, and that, by the uniqueness of $y^{*}$, setting $r_H=1/x_F$ gives
$p(r_H)=F$. Hence $\operatorname{tr}J(p(r_H))=\partial_x f(F)=0$ while
$\det J(p(r_H))>0$, so the eigenvalues are $\pm i\omega$ as claimed.
For item~3, implicit differentiation of $f(1/r,y^{*}(r))\equiv 0$ gives
$\frac{dy^{*}}{dr}=\partial_x f/\left(r^{2}\,\partial_y f\right)$,
which vanishes at $r=r_H$; therefore
\begin{equation*}
\frac{d}{dr}\operatorname{tr}J(p(r))\Big\vert_{r=r_H}
=\partial_{xx}f(F)\cdot\Big(-\frac{1}{r_H^{2}}\Big)
+\partial_{xy}f(F)\cdot\frac{dy^{*}}{dr}\Big\vert_{r=r_H}
=-\frac{\partial_{xx}f(F)}{r_H^{2}}.
\end{equation*}
The claim follows from the standard Hopf bifurcation conditions.
\end{proof}

\begin{figure}[htbp!]
\centering
\begin{tikzpicture}
    \node at (-4,0){
    \begin{tikzpicture}
        \node at (0,0){
    \includegraphics[scale=1.1]{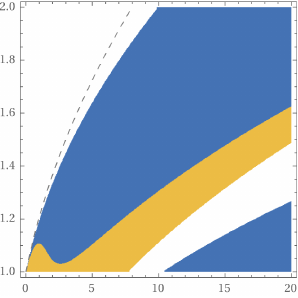}};
    \node at (0,-3){$c$};
    \node at (-3.25,0){$d$};
    \node at (0,3){$\alpha=2\ln3,\;\gamma=\tfrac12$};
    \end{tikzpicture}
    };
    \node at (4,0){
    \begin{tikzpicture}
        \node at (0,0){
    \includegraphics[scale=1.1]{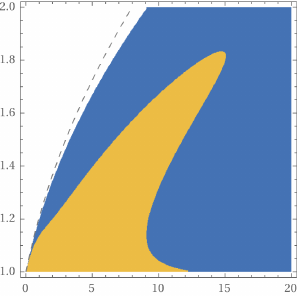}};
    \node at (0,-3){$c$};
    \node at (-3.25,0){$d$};
    \node at (0,3){$\alpha=\tfrac12\ln3,\;\gamma=\tfrac12$};
    \end{tikzpicture}
    };
\end{tikzpicture}
\caption{Criticality of the singular Hopf bifurcation at the leftmost
fold $F_1$ of $\mathcal{C}_0$, in the $(c,d)$-plane at $\gamma=1/2$, for
$\alpha=2\ln 3$ (left) and $\alpha=\tfrac12\ln 3$ (right); compare with the
fold-point diagrams of Figure~\ref{fig:fold-diagrams}. Colour shows the
sign of the intrinsic criticality invariant $\sigma$ of
\cite{de2021intrinsic}, evaluated at $F_1$ with $r=r_H=1/x_{F_1}$
(Proposition~\ref{prop:hopf}): blue, supercritical ($\sigma<0$); yellow,
subcritical ($\sigma>0$). The dashed grey curve is the sufficient
condition $c=4d(d-1)$ of Proposition~\ref{prop:NH_simple}; white marks
parameter values where $\mathcal{C}_0$ has zero (upper left) or four
(lower right) folds, for which the criticality computation is not performed. The boundary along which $\sigma$ changes sign separates supercritical from
subcritical singular Hopf bifurcations; the mechanism driving this change of criticality
is left open, see Section~\ref{Sec:Discussion}.}
\label{fig:criticalitymaps}
\end{figure}

Proposition~\ref{prop:hopf} has the following consequences, which we
highlight individually:
\begin{itemize}
\item  In contrast with the generic
singular Hopf scenario, in which the bifurcation occurs at an
$\mathcal{O}(\varepsilon)$-distance from the parameter value at which
the equilibrium sits exactly on the fold
\cite{Ref:Krupa2001,Ref:Kuehn2015}, here the bifurcation parameter
value carries no $\varepsilon$-correction at all. This is a consequence
of the fact that $g(x,y)=y(1-rx)$ and $\partial_y g$ vanish
simultaneously at the equilibrium.
\item Item~3 connects this analysis to the
degenerate folds studied in Section~\ref{Sec:Results}. At a non-generic
fold ($k> 1$) the transversality of the eigenvalue crossing fails, and the Hopf bifurcation itself degenerates.
\item Since each fold point yields
one Hopf bifurcation, the fold-point diagrams of Figure \ref{fig:fold-diagrams}
also account for the singular Hopf bifurcations of the full system
as $r$ varies.
\item The criticality of these Hopf bifurcations is governed by the
first Lyapunov coefficient, which due to the transcendental nature of
$f$ does not admit a closed form. We therefore evaluate it numerically
through the intrinsic, coordinate-free criticality invariant $\sigma$
of \cite{de2021intrinsic} (with $\sigma<0$ corresponding to a
supercritical bifurcation), which avoids any normal-form transformation
and depends only on Lie derivatives of the fast field at the fold. The
resulting criticality maps are shown in
Figure~\ref{fig:criticalitymaps}. At the parameters used for the
controllers in Section~\ref{Sec:Results} the bifurcation is
supercritical over most of the two-fold region, but its sign changes across an interior curve. The corresponding mechanism by which criticality changes is out of the scope of the current presentation.
\end{itemize}

\begin{figure}[htbp]
    \centering
    \begin{tikzpicture}
        \node at (0,0){
\includegraphics[scale=1.1]{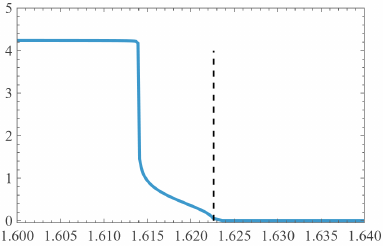}
        };
        \node at (0,-2.5){$r$};
        \node at (-4,0)[rotate=90]{amplitude};
    \end{tikzpicture}
    \caption{Bifurcation diagram showing the singular Hopf bifurcation and the ensuing canard explosion. The parameters used are $\alpha = 2.0$, $\beta = 0.75$, $\gamma = 0.5$, $c = 2.5$, $d = 1.18$, $b = 30.0$, $\varepsilon = 0.001$ with the bifurcation occuring at the left most fold $F_1$ of Figure \ref{Fig:CompetitiveSystemOpenLoop}. The dashed line indicates the bifurcation value $r_H$. Notice how, to the left of the bifurcation value, the amplitude stays small for some interval of $r$, corresponding to canard cycles, and then the amplitude quickly increases (canard explosion).}
    \label{fig:canardexplotion}
\end{figure}

\begin{remark}
    The fold-point diagrams and critical manifolds in Figures \ref{fig:fold-diagrams} and \ref{fig:criticalmanifolds}
describe the singular limit $\varepsilon=0$. For the full system with
$0<\varepsilon\ll 1$, Proposition~\ref{prop:hopf} and \cite{Ref:Krupa2001} (see also Appendix~A) imply
the following: each generic fold point $F$ of $\mathcal{C}_0$ gives
rise to a singular Hopf bifurcation of the full system, located exactly
at $r_H=1/x_F$, at which the equilibrium $(1/r,y^{*}(r))$ crosses the
fold. The ensuing canard explosion produces a family of canard cycles
whose amplitude grows from $\mathcal{O}(\sqrt{\varepsilon})$ to
$\mathcal{O}(1)$ (relaxation oscillation) as $r$ traverses an interval
of exponentially small width $\mathcal{O}(e^{-K/\varepsilon})$, $K>0$,
located at an $\mathcal{O}(\varepsilon)$-distance from $r_H$, see Figure \ref{fig:canardexplotion}. In the uncontrolled system, these canard cycles are therefore
structurally unstable, that is, they exist only for exponentially fine-tuned
values of $r$ and are destroyed by arbitrarily small perturbations. The
fold-point diagrams above identify the regions in parameter space where
such canard cycles can arise. The purpose of the controllers developed
in Section \ref{Sec:Results} is to stabilise a chosen canard cycle $\gamma_h$ as a robust attractor of the closed-loop system.
\end{remark}

\begin{remark}\label{rmk:canardvscycle}
Near a fold point, a standard (relaxation) limit cycle and a canard cycle are
distinguished by their singular limit $\varepsilon\to0$, also known as the \emph{candidate orbit} \cite{Ref:Kuehn2015}, which in this case corresponds to a
closed concatenation of slow arcs on $\mathcal{C}_0$ and fast fibres. A standard
relaxation limit cycle has a candidate that leaves $\mathcal{C}_0$ along a fast
fibre as soon as it reaches the fold, visiting only attracting arcs. A \emph{canard
cycle} is a limit cycle whose candidate contains a \emph{singular canard}, i.e.\ an
arc of the repelling part of $\mathcal{C}_0$.
\end{remark}

\section{Main Results}\label{Sec:Results}

    In this section, we detail the local analysis and the derivation of fast--slow controllers that stabilize canard orbits in a neighborhood of a fold point of the decision--making system \eqref{Eq:OriginalSystem}. We note that equilibrium stabilization near a fold can be viewed as a singular limit of canard stabilization, but the latter provides substantially improved robustness with respect to parameter uncertainty and structural perturbations. Indeed, stabilizing an equilibrium at a fold by tuning a single
parameter (e.g.\ the harvesting rate~$r$) requires precise knowledge of
the fold location: arbitrarily small parameter errors can move the
equilibrium onto the repelling branch of~$\mathcal{C}_0$, immediately
triggering large relaxation oscillations.  By contrast, the canard
controller targets the transition mechanism itself and remains effective
under moderate uncertainty; moreover, it acts on the agents'
decision-making rules rather than on the resource dynamics directly,
making it less intrusive from an applied perspective.
    
    First, we study a general fast--slow system and design compatible controllers capable of stabilizing canards around a not necessarily generic fold point, results summarized in Theorem \ref{Thm:GeneralizedQuadraticControl}, and extending the application of our control methods to more degenerate cases. Although degenerate folds are nongeneric, they naturally arise as effective models when higher-order terms dominate the local geometry near a fold over finite parameter ranges. Moreover, depending on the dynamical features of the system at hand, we provide two different alternatives to accomplish the stabilization of canards, when having either a semi--actuated or a fully controlled system. Thereafter, we employ our general results to system \eqref{Eq:OriginalSystem} in order to stabilize canard orbits in two different fold points, entitling the controlling entity with the capability to stay near the bifurcation or produce abrupt transitions by exploiting the system's criticality. In particular, in our first example, detailed in Section \ref{SubSec:MaxStock}, we stabilize orbits near a fold point that maximizes the stock of renewable resource $y$ by restricting the share of agents exploiting it. On the other hand, our second example, explained in Section \ref{SubSec:OptimalStrategy}, denotes an \emph{optimal exploitation strategy} in which the maximum number possible of agents $x$ benefit from consuming the limited resource $y$ while reducing its stock to the local minimum, producing a favorable scenario both for the population and the controlling entity.

    \subsection{Generalization of the normal form controls}
    \label{Sec:Applications}
    We begin by introducing a general fast--slow system in the form
    \begin{align}
        \begin{split}
            \dot{x} &= F(x,y, \varepsilon),\\
            \dot{y} &= \varepsilon G(x,y, \varepsilon),
            \label{Eq:NormalFormApplicationsGeneralVectorField}
        \end{split}
    \end{align}
    where $0<\varepsilon\ll1$ is the timescale separation parameter, while $F:\mathbb{R}^{m}\times\mathbb{R}^{n}\times\mathbb{R}\xrightarrow{}\mathbb{R}^{m}$ and $G:\mathbb{R}^{m}\times\mathbb{R}^{n}\times\mathbb{R}\xrightarrow{}\mathbb{R}^{n}$ are sufficiently smooth functions such that the associated critical manifold $\mathcal{C}_{0}$ of the singularly perturbed vector field \eqref{Eq:NormalFormApplicationsGeneralVectorField} has at least one point $p=\left( x_{p}, y_{p}
  \right)\in\mathcal{C}_{0}$ as a not necessarily generic fold point, i.e., it satisfies the following conditions: 
\begin{enumerate}
  \item[(i)] $F(x_p,y_p,0)=0$ and $G(x_p,y_p,0)\neq 0$
        \quad(the point lies on $\mathcal C_0$ and the slow flow does not vanish);
  \item[(ii)] $\partial_x F(x_p,y_p,0)=\partial_x^2 F(x_p,y_p,0)=\cdots
        =\partial_x^{\,2k-1} F(x_p,y_p,0)=0$, $k\in\mathbb N$
        \quad(the first $2k-1$ fast derivatives vanish);
  \item[(iii)] $\partial_x^{\,2k} F(x_p,y_p,0)\neq 0$
        \quad(leading fast derivative non-zero, contact order $2k$);
  \item[(iv)] $\partial_y F(x_p,y_p,0)\neq 0$
        \quad(nondegeneracy: $\mathcal C_0$ is a graph over $y$ at the fold).
\end{enumerate}
 
\begin{figure}
            \centering
            \includegraphics[width=0.5\linewidth]{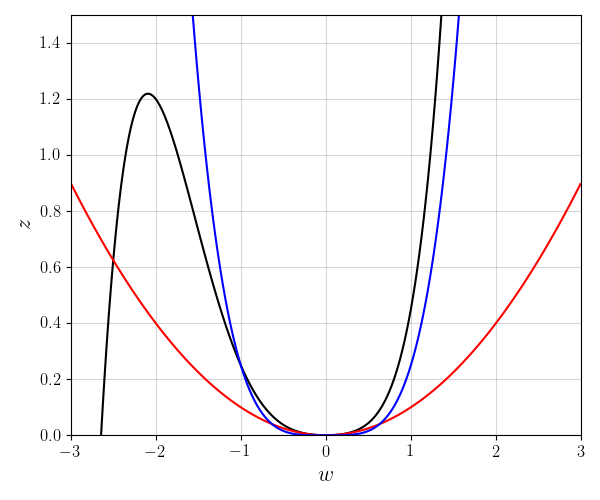}
            \caption{Graphical comparison between two possible approximations to the critical manifold $\mathcal{S}_{0}\coloneqq\left\{ (w,z)\in\mathbb{R}^{2}: 0.1w^{5} + 0.25w^{4} + 0.1w^{2} = z \right\}$ (black) by the manifolds $\mathcal{S}_{2}\coloneqq
            \left\{ (w,z) \in \mathbb{R}^{2}: 0.1w^{2} = z \right\}$ (red) and $\mathcal{S}_{4}\coloneqq
            \left\{ (w,z) \in \mathbb{R}^{2}: 0.25w^{4} = z \right\}$ (blue). As it can be seen, despite that in a small neighborhood of the fold the best approximation is given by $\mathcal{S}_{2}$, as it is the first non--zero derivative, the manifold $\mathcal{S}_{4}$ is indeed a better approximation as the size of the desired canard orbit increases.}
            \label{Fig:DegenerateFold}
        \end{figure}

The parameter $k$ determines the contact order of the parabola--like shape in the
critical manifold of \eqref{Eq:NormalFormApplicationsGeneralVectorField} near the
bifurcation point $p$, which for a generic fold point is $k=1$. It is precisely the
parameter $k$ that enables the use of our control techniques in applications for
which the associated critical manifold $\mathcal{C}_{0}$ is better locally approximated by
a higher--order quadratic--like system ($k>1$) around a not necessarily generic fold point
$p\in\mathcal{C}_{0}$. For instance, in Figure~\ref{Fig:DegenerateFold} the critical
manifold $\mathcal{S}_{0}\coloneqq\left\{ (w,z)\in\mathbb{R}^{2}: 0.1w^{5} + 0.25w^{4} +
0.1w^{2} = z \right\}$ (black) is approximated by $\mathcal{S}_{2}\coloneqq\left\{ (w,z)
\in \mathbb{R}^{2}: 0.1w^{2} = z \right\}$ (red) and $\mathcal{S}_{4}\coloneqq\left\{ (w,z)
\in \mathbb{R}^{2}: 0.25w^{4} = z \right\}$ (blue). The coefficients of
$\mathcal{S}_{0}$ are illustrative and carry no physical meaning; the subscript on each
manifold records the leading power retained, so that $\mathcal{S}_{0}$ is the original,
untruncated manifold, while $\mathcal{S}_{2}$ and $\mathcal{S}_{4}$ are the generic
($k=1$) and quartic ($k=2$) fold approximations. Only even subscripts arise because, by
conditions (ii)--(iii), a fold has even contact order $2k$.
At larger scales $\mathcal{S}_{4}$ provides the better fit, even though $\mathcal{S}_{2}$
represents the first derivative different from zero in the expansion of
\eqref{Eq:NormalFormApplicationsGeneralVectorField} around the fold point located at the
origin. In other words, although a generic approximation provides a better fit in an
arbitrarily small neighbourhood of the fold, higher--order approximations yield a more
accurate representation of the critical manifold over the finite region explored by canard
trajectories. We note that while arbitrary polynomial approximations of the critical
manifold (including linear combinations of different--order truncations) may improve
pointwise accuracy, the control design in Section~\ref{Sec:Results} relies on the
Hamiltonian structure of the leading part, and therefore truncations that destroy this
structure fall outside the present framework.
    \begin{remark}
        The importance of the previous fact in the decision--making model \eqref{Eq:OriginalSystem} is clear as the parameters vary and the critical manifold $\mathcal{C}_{0}$ presents more degenerate forms, especially when there are four fold points. 
        Therefore, depending on each particular case, the critical manifold may be better locally approximated by a higher-order term than by the first derivative different from zero when expanding the fast--slow vector field around the fold point of interest, evidencing the relevance of the contact order parameter $k$ for our results.
    \end{remark}
 
 Hence, as long as conditions \textbf{i}--\textbf{iv} are fulfilled, the dynamics in a region of \eqref{Eq:NormalFormApplicationsGeneralVectorField} sufficiently close to the fold point $p$ are approximated by a system in the form
 \begin{subequations}
 \renewcommand{\theequation}{\theparentequation.\arabic{equation}}
    \begin{align}
            \dot{x} &= a_{c}x^{2k} + b_{c}y + \Tilde{F}\left(x,y,\varepsilon, \alpha(x,y;\lambda)\right), \label{Eq:NormalForm_Fast}\\
            \dot{y} &= -\varepsilon \left( \sigma k a_{c} x^{2k-1} - \alpha(x,y;\lambda) + \Tilde{G}\left(x,y,\varepsilon, \alpha(x,y;\lambda)\right) \right), \label{Eq:NormalForm_Slow}
    \end{align}
    \label{Eq:ApproximationGeneralVectorField}
    \end{subequations}
    with the parabola--like shape and its stability conditions determined by the non--zero constants $a_{c} = \frac{1}{(2k)!}\frac{\partial^{2k} F}{\partial x^{2k}}\left(x_{p}, y_{p}, 0\right)$, and $b_{c} = \frac{\partial F}{\partial y}\left(x_{p}, y_{p}, 0\right)$, obtained through the expansion of \eqref{Eq:NormalFormApplicationsGeneralVectorField} around the fold point $p\in \mathcal{C}_{0}$, and where we reuse the same variables for the coordinate system as in \eqref{Eq:NormalFormApplicationsGeneralVectorField} for simplicity, with the parameter $\sigma = \text{sign}\left(a_{c}b_{c}\right)$. Additionally, the functions $\Tilde{F}$ and $\Tilde{G}$ collect any higher order terms, while the term $\alpha(x,y;\lambda)$ sets the equilibrium point of the slow problem of \eqref{Eq:ApproximationGeneralVectorField}. Furthermore, if we consider that for a neighborhood sufficiently close to the fold point $p$ the higher order terms $\Tilde{F}(x,y,\varepsilon,\alpha(x,y;\lambda))$ and $\Tilde{G}(x,y,\varepsilon,\alpha(x,y;\lambda))$ in \eqref{Eq:ApproximationGeneralVectorField} are zero, then for $\varepsilon>0$ and $\alpha(x,y;\lambda)=0$ orbits of \eqref{Eq:ApproximationGeneralVectorField} are given by level sets of
    \begin{equation}
        H\left( x, y, \varepsilon \right) = \frac{\sigma}{2}e^{2y/\sigma \varepsilon}\left( \frac{b_{c}y}{\varepsilon} + \frac{a_{c}x^{2k}}{\varepsilon} - \frac{\sigma b_{c}}{2}
            \right),
        \label{Eq:HamiltonianGeneralVectorField}
    \end{equation}
   and it is known that canard cycles exist for $H(x,y,\varepsilon)=h$ with $0<-h/b_c<1/4$ \cite{Ref:Krupa2001, Ref:JardonKojakhmetov2022};
the classical range $h\in(0,1/4)$ corresponds to the normalization $(a_c,b_c)=(1,-1)$. Thus, a region of the critical manifold of \eqref{Eq:NormalFormApplicationsGeneralVectorField} sufficiently close to the fold point $p$ can be approximated by a manifold as
    \begin{equation}
        \mathcal{S}_{0} = \left\{ (x,y)\in\mathbb{R}^{2}: a_{c}x^{2k} + b_{c}y + \Tilde{F}(x,y,0,\alpha(x,y;\lambda)) = 0  \right\}.
    \end{equation}
    Now, let us introduce the fast--slow control entries into \eqref{Eq:ApproximationGeneralVectorField}, as
    \begin{align}
            \begin{split}
                \dot{x} &= a_{c}x^{2k} + b_{c}y + \tilde{F}(x, y, \varepsilon, \alpha(x,y;\lambda)) + u(x,y;\xi),\\
                \dot{y} &= -\varepsilon \left( \sigma k a_{c}x^{2k-1} - \alpha(x,y;\lambda) + \tilde{G}(x, y, \varepsilon, \alpha(x,y;\lambda)) + v(x,y;\xi)\right), 
            \end{split}
            \label{Eq:GeneralizedQuadraticSystemHOTControl}
        \end{align}
    where $\xi$ collects all possible controller parameters. The controlled system \eqref{Eq:GeneralizedQuadraticSystemHOTControl} is obtained by adding controllers to each of the equations of the normal form \eqref{Eq:ApproximationGeneralVectorField}; as such, each controller is added to actuate at different scales; $u$ is the fast controller (or acts on the fast variable), while $v$ is in the slow regime. Moreover, observe that without control restrictions both $u(x,y;\xi)$ and $v(x,y;\xi)$ are able to compensate the higher order terms $\Tilde{F}(x,y,\varepsilon,\alpha(x,y;\lambda))$ and $\Tilde{G}(x,y,\varepsilon,\alpha(x,y;\lambda))$, respectively, when enough knowledge on these higher order terms is at hand. Hence, assuming that the higher order terms are compensated by the control entries $u(x,y;\xi)$ and $v(x,y;\xi)$, two possible alternatives arise in order to bring the slow dynamics of  \eqref{Eq:GeneralizedQuadraticSystemHOTControl} to the necessary form depending on the $\alpha(x,y;\lambda)$ term, that is, setting the equilibrium of the slow dynamics in the origin so that a canard point is obtained. First, if $\alpha(x,y;\lambda)=\alpha$ is a constant value, then a coordinate transformation $\hat{x} = x-\alpha$ is enough to set the equilibrium of the slow dynamics of  \eqref{Eq:GeneralizedQuadraticSystemHOTControl} at the origin. On the other hand, if $\alpha(x,y;\lambda)$ presents a different variable or parameter dependence, then the slow control $v(x,y;\xi)$ is able to compensate its effect. Hence, the fast controller $u(x,y;\xi)$ acts as the main control, while the slow input $v(x,y;\xi)$ becomes a support control when the slow dynamics require it. Thus, after compensating the effect of the higher order terms $\Tilde{F}(x,y,\varepsilon,\alpha(x,y;\lambda))$, and $\Tilde{G}(x,y,\varepsilon,\alpha(x,y;\lambda))$, as well as $\alpha(x,y;\lambda)$, the control problem reduces to
    \begin{align}
        \begin{split}
            \dot{x} &= a_{c}x^{2k} + b_{c}y + u(x,y;\xi),  \\
            \dot{y} &= -\varepsilon \sigma k a_{c}x^{2k-1}.
        \end{split}
    \label{Eq:GeneralizedQuadraticSystem}
    \end{align}
    Then, the critical manifold of \eqref{Eq:GeneralizedQuadraticSystem} is $\mathcal{C}_{0} \coloneqq \{ (x, y) \in \mathbb{R}^{2}: y = -(a_{c}/b_{c})x^{2k} \}$, with stability conditions given by the sign of $a_{c}$ and $u(x,y;\xi)$ being a compatible controller. Moreover, notice that \eqref{Eq:GeneralizedQuadraticSystem} defines a conservative system with the Hamiltonian \eqref{Eq:HamiltonianGeneralVectorField}.
    
    Now, the aim of the control strategy is to stabilize one orbit of \eqref{Eq:GeneralizedQuadraticSystem}, namely $\gamma_{h}=\{ (x,y) \in \mathbb{R}^{2} : H=h \}$, related to the canard trajectory of interest. Hence, the control objective is to reduce the error existing between a level set of \eqref{Eq:HamiltonianGeneralVectorField} and the desired trajectory $\gamma_{h}$, determined by the constant $h$. Therefore, we want to stabilize the error dynamics $\Tilde{H} \coloneqq H - h$ to zero. Specifically, we select $h = -\tfrac{1}{4}\,\mathrm{sign}(b_c)\,e^{-c_c/\varepsilon}$ with $c_c \in (0,\infty)$.
Then $-h/b_c = \tfrac{1}{4|b_c|}\,e^{-c_c/\varepsilon} > 0$, and $-h/b_c < 1/4$ whenever
$e^{-c_c/\varepsilon} < |b_c|$, so for $\varepsilon$ small enough $h$ lies in the canard-cycle
range $0 < -h/b_c < 1/4$. Note that this selection of $h$ is done in order to improve the precision by using exponential variations. Thus, it is clear that
    \begin{equation}
            \dot{\Tilde{H}} = \frac{e^{2y/\sigma \varepsilon}}{\varepsilon^{2}}\left( \dot{y}\left( a_{c}x^{2k} + b_{c}y \right) + \varepsilon \sigma k a_{c}x^{2k-1}\dot{x} \right).
            \label{Eq:HDerivativeRaw}
        \end{equation}
        Then, by substituting \eqref{Eq:GeneralizedQuadraticSystem} in \eqref{Eq:HDerivativeRaw} yields
        \begin{equation}
            \dot{\tilde{H}} = \frac{e^{2y/\sigma \varepsilon}}{\varepsilon}\left(\sigma k a_{c}x^{2k-1}u \right).
            \label{Eq:HDerivativeControl}
        \end{equation}
        Subsequently, an alternative to guarantee the local asymptotic stability of the desired level set $\gamma_{h}$ is to define the error dynamics as $\dot{\Tilde{H}} = -A_{c}\Tilde{H} = -A_{c}(H - h)$, for some $A_{c} > 0$. Thus, the fast control entry is given by
        \begin{equation}
            u(x,y;\xi) = -\frac{\varepsilon B_{c} x}{\sigma k a_{c}}(H-h)e^{-2y/\sigma \varepsilon},
            \label{Eq:FastControlEntry}
        \end{equation}
        where we have desingularized the origin by setting $A_{c} = B_{c}x^{2k}$, positive for every $x\in\mathbb{R}\backslash\{0\}$, and $B_{c}>0$ is the controllers' gain. 
        
\begin{figure*}[htbp]
            \begin{subfigure}[htbp]{0.5\textwidth} \includegraphics[width=1.0\linewidth]{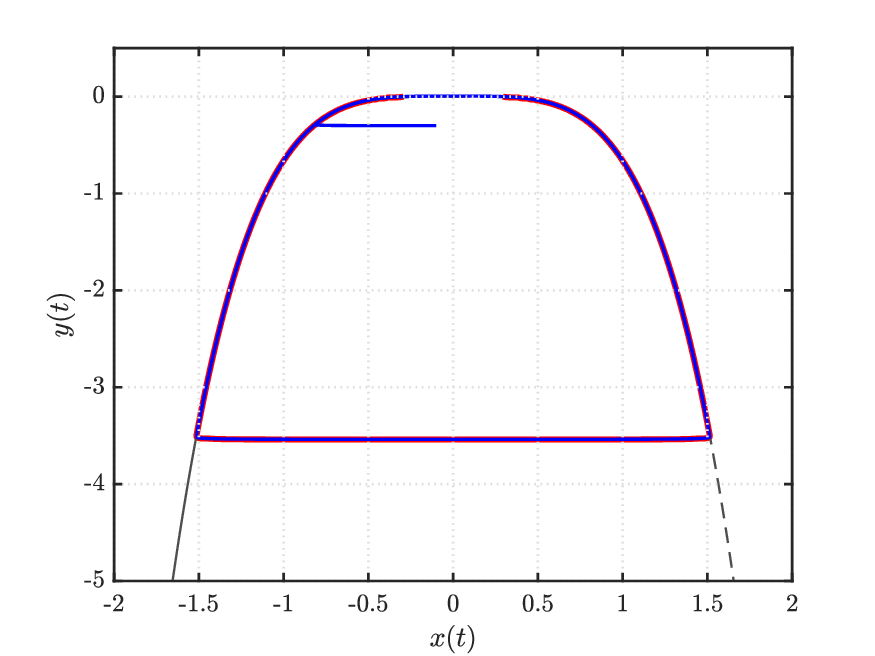}
            \end{subfigure}\hfill
            \begin{subfigure}[htbp]{0.5\textwidth}
            \includegraphics[width=1.0\linewidth]{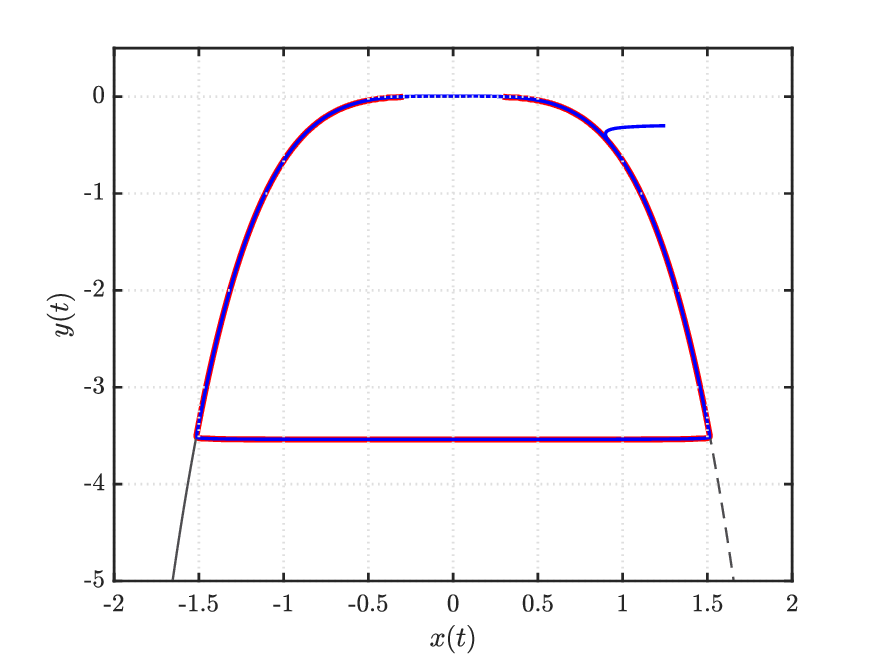}
            \end{subfigure}
            \vskip \baselineskip
            \begin{subfigure}[htbp]{0.5\textwidth}
            \includegraphics[width=1.0\linewidth]{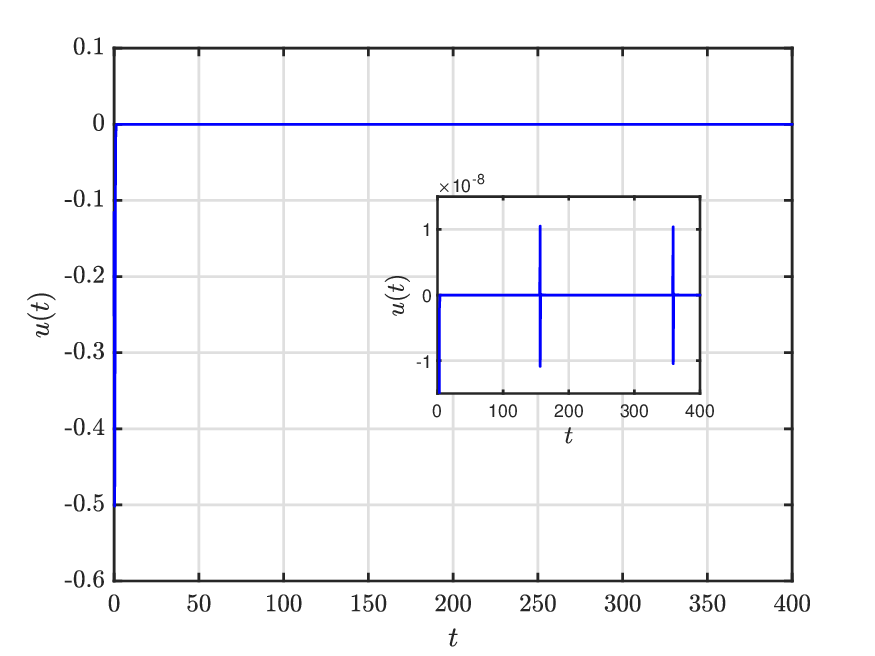}
            \end{subfigure}\hfill
            \begin{subfigure}[htbp]{0.5\textwidth}
            \includegraphics[width=1.0\linewidth]{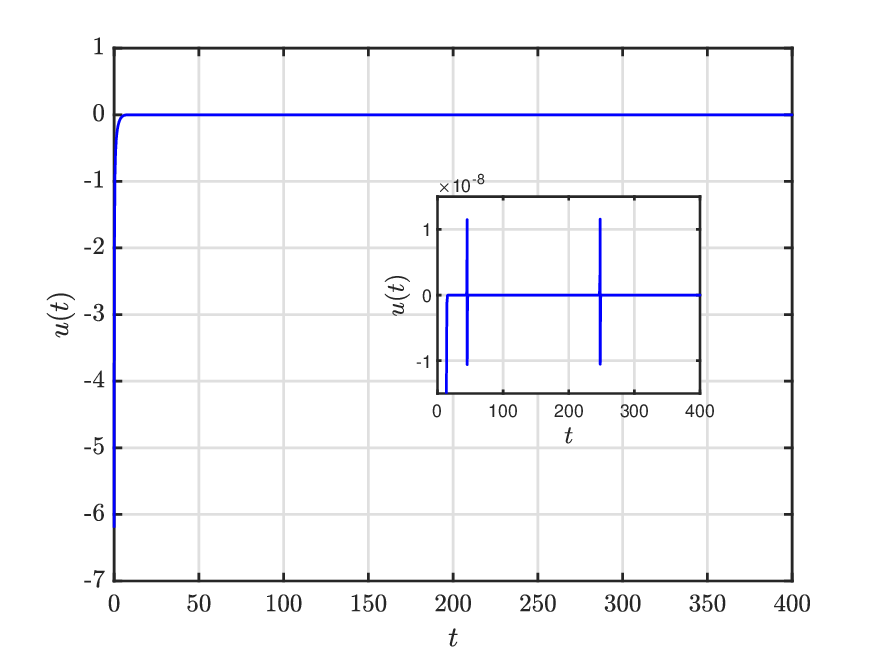}
            \end{subfigure}
            \vskip \baselineskip
            \begin{subfigure}[htbp]{0.5\textwidth}
            \includegraphics[width=1.0\linewidth]{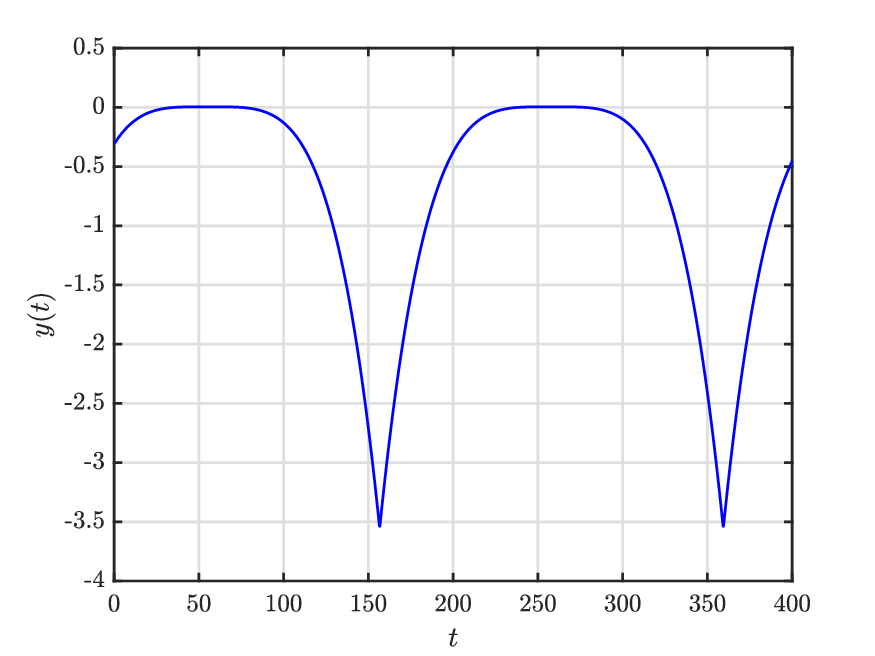}
            \end{subfigure}\hfill
            \begin{subfigure}[htbp]{0.5\textwidth}
            \includegraphics[width=1.0\linewidth]{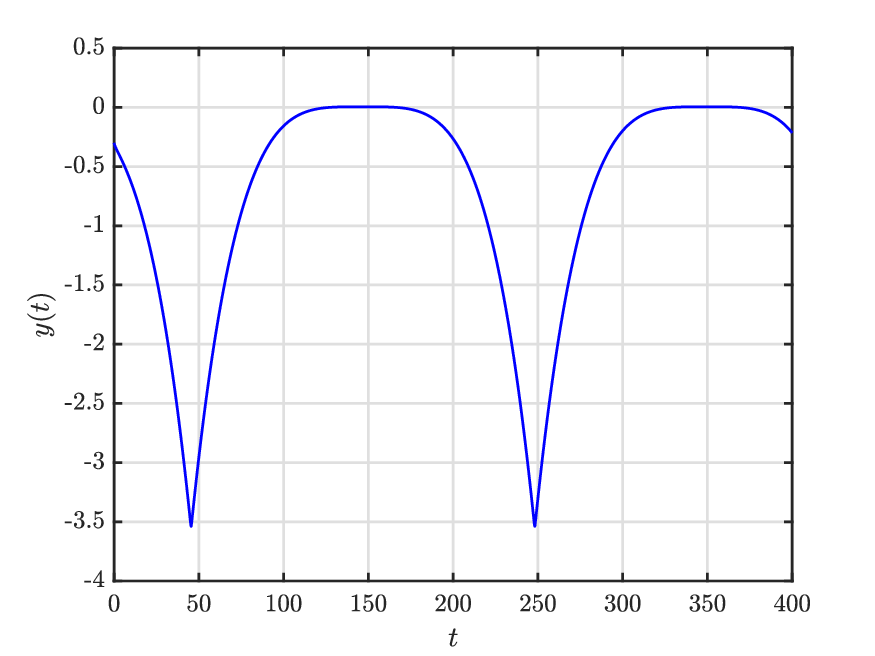}
            \end{subfigure}
            \caption{$x-y$ plane (top), control output (center), and $y$--time series (bottom) for a trajectory of \eqref{Eq:GeneralizedQuadraticSystem} (solid blue), with the fast controller  \eqref{Eq:FastControlEntry}, for initial conditions $(x(0), y(0)) = (-0.1, -0.3)$, and $(x(0), y(0)) = (1.25, -0.3)$, inside (left) and outside (right) of the target level curve, respectively.
        For the $x-y$ plane, the attracting and repelling branches of the critical manifold $\mathcal{C}_{0}$ are depicted in solid and dashed black, respectively, while the target level set $\gamma_{h}$ is presented in solid red. The center panel in the control input shows that $u(t)$ is negligible away from the fold and fires only in brief bursts as the orbit rounds it (magnified inset).} Parameters: $a_{c} = 2.0$, $b_{c} = 3.0$, $c_{c} = 7.0$, $B_{c} = 10.0$, $\varepsilon = 0.01$, and $k = 2$.
        \label{Fig:FastControlQuadraticSystem}
        \end{figure*}

        To show stability, we define a candidate Lyapunov function as
        \begin{equation}
            L(x, y) = \frac{1}{2}\Tilde{H}^{2}.
            \label{Eq:LyapunovGeneralControl}
        \end{equation}
            Observe that \eqref{Eq:LyapunovGeneralControl} is positive for every $\Tilde{H}\neq 0$, and that $L=0$ if and only if $\Tilde{H}=0$, if and only if $(x, y)\in \gamma_{h}$, which we recall is the control target defined earlier. Then, it is easy to check that
        \begin{equation}
            \dot{L} = \Tilde{H} \dot{\Tilde{H}} = -B_{c}x^{2k}\Tilde{H}^{2}.
            \label{Eq:LyapunovFastControlDerivative}
        \end{equation}
        Finally, to demonstrate the asymptotic stability of $\gamma_{h}$ as \eqref{Eq:LyapunovFastControlDerivative} is only negative semidefinite, by LaSalle's invariance principle \cite{Ref:Lasalle1960}, trajectories of the fast--slow vector field \eqref{Eq:GeneralizedQuadraticSystem} under the control action $u(x,y;\xi)$ as \eqref{Eq:FastControlEntry} reach, in finite time, the largest invariant set contained in
        \begin{equation}
            \mathcal{I} = \left\{ (x, y)\in \mathbb{R}^{2}:\dot{L} = 0 \right\} = \{ x=0 \} \cup \left\{ \Tilde{H} = 0 \right\}.
            \label{Eq:LyapunovInvariantSet}
        \end{equation}
        It is important to mention that, as long as $y\neq0$, the vector field \eqref{Eq:GeneralizedQuadraticSystem} does not vanish when it reaches \eqref{Eq:LyapunovInvariantSet}, allowing the usage of the fast control \eqref{Eq:FastControlEntry} in the quadratic--like system \eqref{Eq:GeneralizedQuadraticSystem}. Notice, however, that ${x=0}$ is generically not invariant for the closed--loop dynamics \eqref{Eq:GeneralizedQuadraticSystem}. Actually, setting $x=0$ in \eqref{Eq:GeneralizedQuadraticSystem} reduces to $\left( \dot{x}, \dot{y} \right) = \left( b_{c}y, 0 \right)$, inducing an increment or decrement in the share of agents $x$ depending of the values of $b_{c}$. Hence, trajectories of \eqref{Eq:GeneralizedQuadraticSystem} eventually reach $\mathcal{I} = \{ (x,y) = (0,0)\} \cup \{ \Tilde{H} = 0 \}$, and since the origin is an unstable equilibrium point of \eqref{Eq:GeneralizedQuadraticSystem} with the controller as \eqref{Eq:FastControlEntry}, we have that every trajectory with initial condition different from zero eventually reaches the set $\{ \Tilde{H} = 0\}$ as $t\xrightarrow{}\infty$. In Figure \ref{Fig:FastControlQuadraticSystem} we present an example of the effect that the fast control entry \eqref{Eq:FastControlEntry} has in system \eqref{Eq:GeneralizedQuadraticSystem}. The target level set $\gamma_{h}$ and the resulting trajectory are shown in solid red and blue, respectively. As it can be appreciated, the fast control scheme makes the system to behave in the desired manner, causing the orbit to follow unstable branches of the critical manifold $\mathcal{C}_{0}$ for time of order $\mathcal{O}(1)$. Hence, canard orbits originating from a fold point in system \eqref{Eq:GeneralizedQuadraticSystem} are stable when implementing an adequate fast control law as \eqref{Eq:FastControlEntry}. 
        
        \begin{remark}
            The canard cycles we wish to control are the level sets $\{H=h\}$. Theorem \ref{Thm:GeneralizedQuadraticControl} turns a chosen reference cycle into a robust attractor. Particularly, it
provides a compatible feedback $u$ that makes the target level
set $\gamma_h$ locally asymptotically stable near the fold, including non-generic folds
($k>1$), see Figure \ref{Fig:FastControlQuadraticSystem}. 
        \end{remark}

        We summarize these results in the following theorem.

    \begin{theorem}[Generalized quadratic system control]
    \label{Thm:GeneralizedQuadraticControl} Let a general fast--slow quadratic--like system \eqref{Eq:NormalFormApplicationsGeneralVectorField} be given locally as \eqref{Eq:ApproximationGeneralVectorField}, with an associated critical manifold $\mathcal{C}_{0}$, and let the Hamiltonian $H(x,y,\varepsilon)$ be defined by \eqref{Eq:HamiltonianGeneralVectorField}. Then, if the higher order terms in the expansion \eqref{Eq:ApproximationGeneralVectorField} are sufficiently small, the compatible controller $u(x,y;\xi) = -\frac{\varepsilon B_{c} x}{\sigma k a_{c}}(H-h)\exp\left({-2y/\sigma \varepsilon}\right)$ renders the target orbit $\gamma_{h} = \{ (x, y)\in\mathbb{R}^{2}:H=h \}$ locally asymptotically stable, stabilizing canard orbits in a vicinity of the origin of \eqref{Eq:ApproximationGeneralVectorField}, 
    with $a_{c}=\frac{1}{(2k)!}\frac{\partial^{2k}F}{\partial x^{2k}}(x_p,y_p,0)$, and $b_{c} = \frac{\partial F}{\partial y}(x_p,y_p,0)$, for $h$ with $0<-h/b_c<1/4$, $k\in\mathbb{N}$, $B_{c}>0$, $0<\varepsilon \ll 1$, and $\sigma = \emph{sign}(a_{c}b_{c})$. Moreover, under the appropriate translations $x\mapsto x+x_p$, $y\mapsto y+y_p$ this controller stabilizes an equivalent canard near $(x_p,y_p)$ of \eqref{Eq:NormalFormApplicationsGeneralVectorField}.
    \end{theorem}
    The proof of Theorem \ref{Thm:GeneralizedQuadraticControl} directly follows the previously presented analysis and, as long as the higher order terms $\Tilde{F}(x,y,\varepsilon,\alpha(x,y;\lambda))$ and $\Tilde{G}(x,y,\varepsilon,\alpha(x,y;\lambda))$ of \eqref{Eq:ApproximationGeneralVectorField} are sufficiently small, the control scheme synthesized in Theorem \ref{Thm:GeneralizedQuadraticControl} produces canard cycles in a neighbourhood of the fold point $p=(x_{p}, y_{p})\in\mathcal{C}_{0}$ of \eqref{Eq:NormalFormApplicationsGeneralVectorField}. This fact shows the true strength of our control technique since, regardless of the complexity of system \eqref{Eq:NormalFormApplicationsGeneralVectorField}, as long as its critical manifold satisfies the previously stated conditions, it is possible to stabilize canard orbits in a neighborhood of the fold point $p$, even if \eqref{Eq:NormalFormApplicationsGeneralVectorField} is not explicitly given in normal form. Finally, it is noteworthy to mention that as the amplitude of the target level set $\gamma_{h} = \{ (x,y)\in\mathbb{R}^{2}:H=h \}$ increases, the approximation error between the associated critical manifold $\mathcal{C}_{0}$ and the target level set $\gamma_{h}$ will accordingly increase, producing an undesired response in \eqref{Eq:NormalFormApplicationsGeneralVectorField} caused by the disparity of both trajectories due to the effect of the higher order terms. To improve the aforementioned behaviour, a possible alternative is the use of complementary control schemes acting sufficiently far away from the bifurcation such that when the trajectory travels in a vicinity that is correctly approximated by our expansion the system is under the effect of our control scheme, but once the expansion error increases beyond a tolerance threshold, an additional controller redirect the solution in order to follow the flow along the critical manifold $\mathcal{C}_{0}$.

Next, we discuss on the effectiveness of the fast--slow control technique developed in this section for the control of canard orbits near a fold point in the presence of parametric uncertainties. 

\subsection{Control of canards in the perturbed normal form}
    \label{SubSec:PerturbationNormalForm}
\begin{figure*}[htbp]
            \begin{subfigure}[htbp]{0.49\textwidth}
                \includegraphics[width=0.97\linewidth]{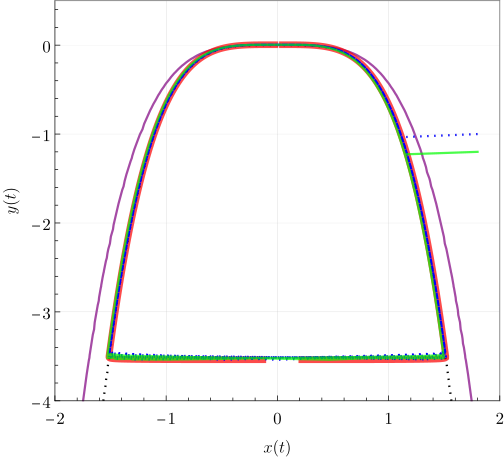}
                \end{subfigure}\hfill
                \begin{subfigure}[htbp]{0.49\textwidth}
                \includegraphics[width=1\linewidth]{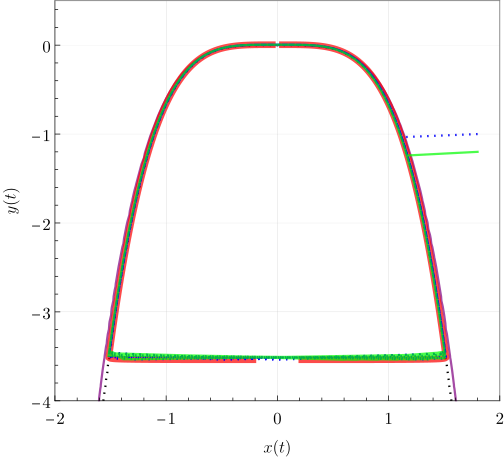}
                \end{subfigure}
                \vskip \baselineskip
                \begin{subfigure}[htbp]{0.49\textwidth}
                \includegraphics[width=0.97\linewidth]{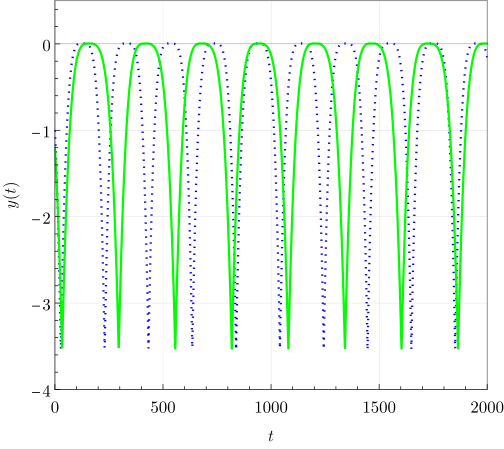}
                \end{subfigure}\hfill
                \begin{subfigure}[htbp]{0.49\textwidth}
                \includegraphics[width=1\linewidth]{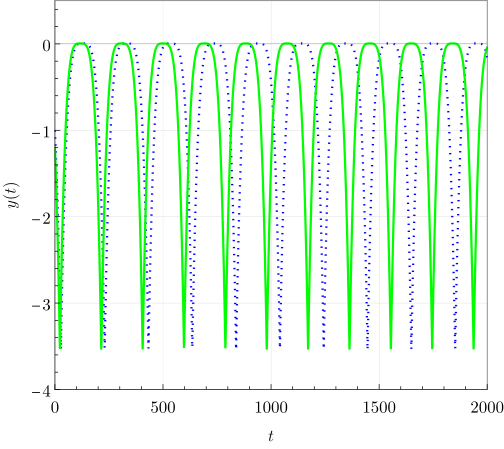}
                \end{subfigure}
                \caption{
                Stabilization of canard cycles in the perturbed normal form \eqref{Eq:PerturbedNormalForm} using the controller \eqref{Eq:FastControlEntry} designed on the unperturbed form \eqref{Eq:GeneralizedQuadraticSystem}, for $\delta_{b}=0.5$ with $\delta_{a}=-0.5$ (left) and $\delta_{a}=0.1$ (right), both satisfying \eqref{Eq:ConditionPerturbationStability}. Upper row (phase portraits): unperturbed manifold and response in dotted black/blue, perturbed in solid purple/green; target level set $\gamma_{h}=\{(x,y)\in\mathbb{R}^{2}:H=h\}$ of the unperturbed Hamiltonian \eqref{Eq:HamiltonianGeneralVectorField} in red. Although \eqref{Eq:FastControlEntry} uses only unperturbed information, it stabilizes the canard cycle in both cases; the lower row ($y(t)$) shows only a phase shift, i.e.\ a speed change that does not alter the stabilization. Parameters: $a=2.0$, $b=3.0$, $\delta_{b}=0.5$, $c=7.0$, $B=50.0$, $h=-2.46492\times10^{-305}$, $k=2$, $\varepsilon=0.01$; $\delta_{a}=-0.5$ (left), $\delta_{a}=0.1$ (right).}
                
                \label{Fig:PerturbationsNormalForm}
            \end{figure*}

    Consider the control problem \eqref{Eq:GeneralizedQuadraticSystem} under the effect of parametric perturbations in the form
    \begin{align}
        \begin{split}
            \dot{x} &= \left(a_{c} + \delta_{a}\right)x^{2k} + \left( b_{c} + \delta_{b} \right)y + u(x,y;\xi),\\
            \dot{y} &= -\varepsilon \sigma k\left(a_{c} + \delta_{a}\right)x^{2k-1},
            \label{Eq:PerturbedNormalForm}
        \end{split}
    \end{align}
    where $\delta_{a}$ and $\delta_{b}$ represent sufficiently small modelling uncertainties in \eqref{Eq:GeneralizedQuadraticSystem}, with the same parameter interpretation as before. Hence, orbits of \eqref{Eq:PerturbedNormalForm} are given by level sets of
    \begin{equation}
        H_{p}\left( x, y, \varepsilon \right) = \frac{\sigma}{2}e^{2y/\sigma \varepsilon}\left( \frac{\left(b_{c}+\delta_{b}\right)}{\varepsilon}y + \frac{\left(a_{c}+\delta_{a}\right)}{\varepsilon}x^{2k} - \frac{\sigma \left(b_{c}+\delta_{b}\right)}{2}
            \right),
        \label{Eq:HamiltonianPerturbedNormalForm}
    \end{equation}
    with the subscript $p$ standing for \emph{perturbation} in $H_{p}(x,y,\varepsilon)$. Therefore, our control objective is to stabilize canard orbits given by level sets of \eqref{Eq:HamiltonianPerturbedNormalForm} in a vicinity of the fold point at the origin of \eqref{Eq:PerturbedNormalForm}, by implementing the compatible fast--slow controller \eqref{Eq:FastControlEntry}, designed in the unperturbed normal form \eqref{Eq:GeneralizedQuadraticSystem} and without explicit information regarding the perturbations $\delta_{a}$ and $\delta_{b}$. We begin by defining the error $\Tilde{H}_{p} = H_{p} - h$, with the associated perturbed error dynamics given by
    \begin{equation}
        \dot{\Tilde{H}}_{p} = \frac{e^{2y/\sigma \varepsilon}}{\varepsilon} \sigma k \left( a_{c} + \delta_{a} \right)x^{2k-1}u, 
        \label{Eq:HDerivativePerturbed}
    \end{equation}
    where we follow the same procedure as for the unperturbed case, aiming to stabilize the error dynamics $\dot{\Tilde{H}}_{p}$ to zero. By substituting the fast--slow controller \eqref{Eq:FastControlEntry} in \eqref{Eq:HDerivativePerturbed} yields
    \begin{equation}
        \dot{\Tilde{H}}_{p} = -\frac{a_{c}+\delta_{a}}{a_{c}}B_{c}x^{2k}(H-h),
        \label{Eq:HDerivativePerturbedControl}
    \end{equation}
    with $B_{c}>0$ the controller gain, while $H(x,y,\varepsilon)$ is \eqref{Eq:HamiltonianGeneralVectorField}. Moreover, since the perturbations $\delta_{a}$ and $\delta_{b}$ are assumed to be sufficiently small, we rewrite \eqref{Eq:HDerivativePerturbedControl} as
    \begin{equation}
        \dot{\Tilde{H}}_{p} = -\Tilde{B}_{c}x^{2k}\left( \Tilde{H}_{p} - H_{\delta}\right),
        \label{Eq:HDErivativePerturbedControl2}
    \end{equation}
    where $\Tilde{B}_{c}>0$ is the new controller gain, and
    \begin{equation}
        H_{\delta}(x,y,\varepsilon) = \frac{\sigma}{2}e^{2y/\sigma \varepsilon}\left( \frac{\delta_{b}}{\varepsilon}y + \frac{\delta_{a}}{\varepsilon}x^{2k} - \frac{\sigma \delta_{b}}{2} \right),
        \label{Eq:HamiltonianPerturbation}
    \end{equation}
    since the unperturbed Hamiltonian \eqref{Eq:HamiltonianGeneralVectorField} is equal to the difference between the perturbed Hamiltonian \eqref{Eq:HamiltonianPerturbedNormalForm} and \eqref{Eq:HamiltonianPerturbation}, i.e. $H(x,y,\varepsilon) = H_{p}(x,y,\varepsilon) - H_{\delta}(x,y,\varepsilon)$. Observe from \eqref{Eq:HDErivativePerturbedControl2} that the error dynamics $\dot{\Tilde{H}}_{p} = -\Tilde{B}_{c}x^{2k}\Tilde{H}_{p} + \Tilde{B}_{c}x^{2k}H_{\delta}$ are bounded from above by $-\Tilde{B}_{c}x^{2k}\Tilde{H}_{p}$. Therefore, in order for the perturbed error to converge to zero in finite time, it is sufficient to show that \eqref{Eq:HamiltonianPerturbation} is negative, and restricted to the perturbed critical manifold $y = -\sigma|(a_{c}+\delta_{a})/(b_{c}+\delta_{b})|x^{2k}$, the aforementioned condition is satisfied as long as
    \begin{equation}
        \sigma \delta_{a} < \delta_{b}\left| \frac{a_{c}}{b_{c}} \right|,
    \label{Eq:ConditionPerturbationStability}
    \end{equation}
    for $\delta_{b}>0$. In Figure \ref{Fig:PerturbationsNormalForm}, we show examples for the stabilization of canard cycles in the perturbed normal form \eqref{Eq:PerturbedNormalForm} by implementing the fast--slow controller \eqref{Eq:FastControlEntry} obtained from the unperturbed normal form \eqref{Eq:GeneralizedQuadraticSystem}. Considering $\delta_{b}$ positive, in the left and right columns we present the cases for which $\delta_{a}$ is negative and positive, respectively, such that condition \eqref{Eq:ConditionPerturbationStability} is satisfied. In the upper row the phase portraits are depicted, where the critical manifold and response associated to the unperturbed problem \eqref{Eq:GeneralizedQuadraticSystem} are represented in dotted black and blue, while the critical manifold and solution associated to the perturbed case \eqref{Eq:PerturbedNormalForm} are shown in solid purple and green, correspondingly. Moreover, the target level set $\gamma_{h}=\{ (x,y)\in \mathbb{R}^{2}:H=h \}$ is again presented in solid red, and observe that the stabilization of the canard trajectory is achieved for both cases. Notice that, under condition \eqref{Eq:ConditionPerturbationStability}, the perturbations $\delta_{a}$ and $\delta_{b}$ cause a phase shift between the responses of \eqref{Eq:GeneralizedQuadraticSystem} and \eqref{Eq:PerturbedNormalForm}, effect that can be appreciated in the comparison of the slow dynamics for both systems presented in the lower row of Figure \ref{Fig:PerturbationsNormalForm}. However, such difference does not modify the stabilization of the desired periodic pattern  in the perturbed normal form \eqref{Eq:PerturbedNormalForm}. 

    In what follows, we employ the results summarized in Theorem \ref{Thm:GeneralizedQuadraticControl} to stabilize canard cycles in a neighborhood of a fold point of the main model \eqref{Eq:OriginalSystem}, even in the presence of modelling uncertainties when condition \eqref{Eq:ConditionPerturbationStability} is satisfied.
    
    \subsection{Control of canards in the decision--making system}
    \label{SubSec:MaxStock}

    Let us recall that under assumption A1-A3, the original system \eqref{Eq:OriginalSystem} is reduced to the simpler one \eqref{eq:main2}. For convenience let us recall that \eqref{eq:main2} reads as
    \begin{subequations}
    \renewcommand{\theequation}{\theparentequation.\arabic{equation}}
    \begin{align}
            \dot{x} &= \gamma(1-x)\left( \frac{1}{1+e^{-(\alpha + \delta(x,y))}} \right) - x\left( \frac{1}{1+e^{-(\alpha-\delta(x,y))}} \right),\label{Eq:SimplifiedSystem_Fast}\\
            \dot{y} &= \varepsilon y (1-rx),\label{Eq:SimplifiedSystem_Slow}
    \end{align}
    \label{Eq:SimplifiedSystem}
    \end{subequations}

\begin{figure*}[htbp]
            \begin{subfigure}[htbp]{0.5\textwidth}
                \includegraphics[width=0.9\linewidth]{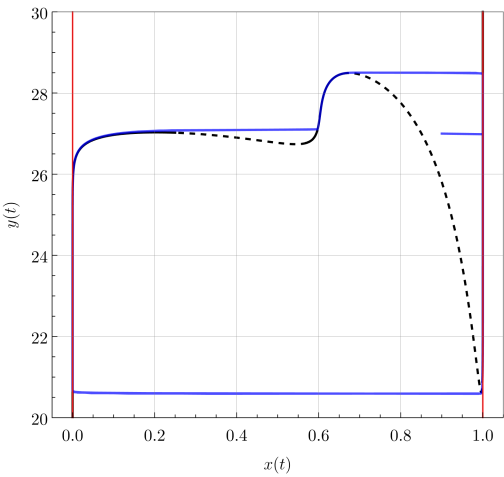}
                \end{subfigure}\hfill
                \begin{subfigure}[htbp]{0.5\textwidth}
                \includegraphics[width=0.9\linewidth]{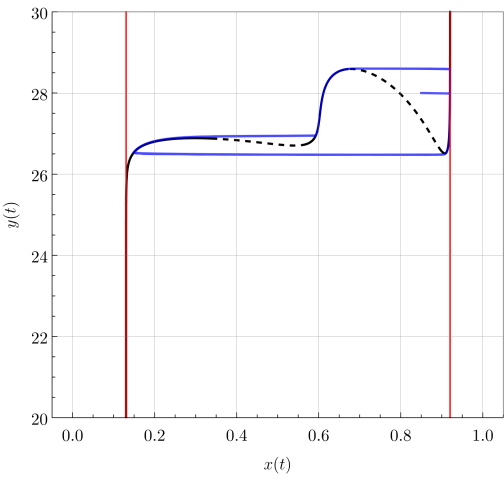}
                \end{subfigure}
                \vskip \baselineskip
                \begin{subfigure}[htbp]{0.5\textwidth}
                \includegraphics[width=0.9\linewidth]{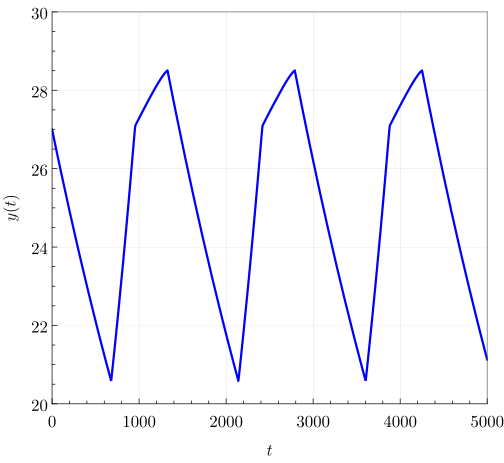}
                \caption{}
                \end{subfigure}\hfill
                \begin{subfigure}[htbp]{0.5\textwidth}
                \includegraphics[width=0.9\linewidth]{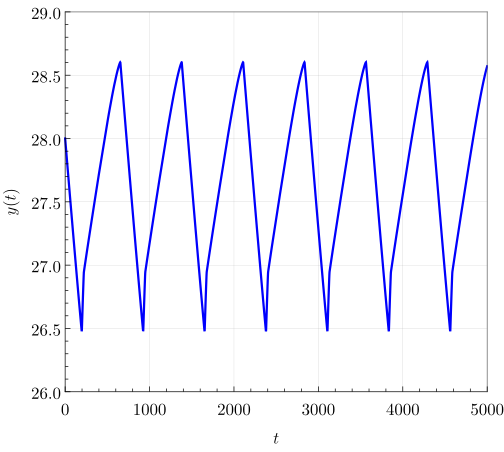}
                \caption{}
                \end{subfigure}
                \caption{Effect of the unconditional exploration values $\eta_{1,2}=0$ (left), and $\eta_{1,2}\in(0,1)$ (right) on the decision--making model \eqref{Eq:OriginalSystem}. Upper row: $x-y$ planes, where stable and unstable regions of the critical manifold are depicted in solid and dashed black, respectively, while the asymptotes and response of \eqref{Eq:OriginalSystem} are shown in red and blue, correspondingly. Lower row: Variations on the slowly renewable resource $y$. Parameters: $\alpha_{1}=0.9$, $\alpha_{2}=2.5$, $\beta_{1}=4.0$, $\beta_{2}=2.0$, $\gamma_{1}=3.0$, $\gamma_{2}=2.0$, $b=30$, $c=1.5$, $d=1.1$, $\varepsilon = 0.001$, and $r=1.4$, with $\eta_{1}=\eta_{2}=0$ (left), and $\eta_{1}=0.1$, $\eta_{2}=0.13$ (right).}
                \label{Fig:AsymptotesOriginalSystem}
            \end{figure*}

    In particular, the choice $\eta_{1}=\eta_{2}=0$ (Assumption A1) allows for the interesting dynamics of \eqref{Eq:OriginalSystem} to occur in the whole domain $x\in[0,1]$. Nevertheless, we emphasize that every pattern observed for $\eta_{1} = \eta_{2} = 0$ is also produced for any $\eta_{1}, \eta_{2} \in (0,1)$ but in a reduced region of the domain of $x$. By performing an asymptotic analysis of system \eqref{Eq:OriginalSystem}, we find the location of the left and right asymptotes of the associated critical manifold $\mathcal{C}_{0}$, being $x_{L} = (\gamma_1 \eta_1)/(\gamma_{1}\eta_{1} + \gamma_2)$, and $x_{R} = \gamma_1/(\gamma_{1}+\gamma_{2}\eta_{2})$, when $y\rightarrow -\infty$ and $y\rightarrow \infty$, respectively. Hence, when $\eta_{1} = \eta_{2} = 1$, the dynamics of \eqref{Eq:OriginalSystem} occur along a vertical line centered at $x = \gamma_{1}/(\gamma_{1} + \gamma_{2})$, dramatically reducing the possible behaviours to be observed. In Figure \ref{Fig:AsymptotesOriginalSystem} we show the effect of varying the unconditional exploration rates $\eta_{1,2}$ when fixing the rest of the parameters in system \eqref{Eq:OriginalSystem}. On the left panel, we consider $\eta_1 = \eta_2 = 0$, which represents a scenario where, once the agents have decided to switch from one strategy to another with rates $\gamma_{1,2}>0$, they change in a completely interested manner with respect to the cost--benefit determined by \eqref{Eq:ProfitDifference}. In contrast, in the right panel we show the case when the agents modify their strategy following a more informed conviction, i.e., $\eta_{1,2}\in(0,1)$. Notice that, by increasing the unconditional exploration rates, the presence of abrupt changes in the share of agents following an exploitation policy is considerably reduced without the need of an external controller. However, although desirable, this scenario is rather unreliable as it completely depends on the agents' goodwill, motivating further the need of an outer controlling entity. Lastly, when every agent switch in an entirely informed way a complete balance between the two groups of agents exploiting each strategy is reached, allowing for the renewable resource only to increase, which would be the case when $\eta_{1,2}=1$, and the share of agents exploiting the limited resource $y$ would be constant with value $x=\gamma_{1}/(\gamma_{1}+\gamma_{2})$.

    In addition, we have identified in Section \ref{Sec:Assumptions} that the critical manifold of \eqref{eq:main2} (and hence of \eqref{Eq:OriginalSystem}) can have $0$, $2$, or $4$ fold points. This motivates us to use the control approach of Theorem  \ref{Thm:GeneralizedQuadraticControl} to stabilize canards close to the fold points. For practical purposes, we will focus on cases where the critical manifold has exactly two fold points.

Therefore, by introducing control inputs in system \eqref{Eq:SimplifiedSystem} yields
    \begin{align}
        \begin{split}
            \dot{x} &= \gamma(1-x)\left( \frac{1}{1+e^{-(\alpha + \delta(x,y))}} \right) - x\left( \frac{1}{1+e^{-(\alpha-\delta(x,y))}} \right) + u(x,y;\xi),\\
            \dot{y} &= \varepsilon\left( y (1-rx) +  v(x,y;\xi) \right),
            \label{Eq:SystemControl}
        \end{split}
    \end{align}
    where $u(x, y;\xi)$ and $v(x,y;\xi)$ (with $v(x,0;\xi)=0$, see Remark \ref{Rem:GraphC0} below) represent the fast and slow control components, respectively, while $\xi$ collects the possible control parameters involved. 
    
    \begin{remark}
    \label{Rem:GraphC0}
        Let us justify the reason for selecting the slow control as mentioned above. Notice that the $y$--dynamics for the open--loop of \eqref{Eq:SimplifiedSystem} is in the form $\dot y=\varepsilon y(1-rx)$. This is not, at first sight, compatible with the normal form \eqref{Eq:ApproximationGeneralVectorField}. However, the reduction to the critical manifold yields slow reduced systems of the form
        \begin{equation}
            y'=g(y) \qquad \text{and} \qquad y'=y(g(y)),
        \end{equation}
        related to \eqref{Eq:NormalForm_Slow}, and  \eqref{Eq:SimplifiedSystem_Slow} respectively. These systems are $C^\infty$--equivalent for $y>0$. To ensure this equivalence, we make sure that the parameters are chosen so that the point $F_2$ (see Figure \ref{Fig:CompetitiveSystemOpenLoop}) is uniformly bounded away from $\{y=0\}$.
    \end{remark}
    
    Note that the selection of the controllers $u(x,y;\xi)$ and $ v(x,y;\xi)$ in \eqref{Eq:SystemControl} is inspired from different physical applications \cite{Ref:Durham2008, Ref:Izhikevich2006, Ref:Ermentrout2010, Ref:JardonKojakhmetov2019}. Particularly, for our decision--making model \eqref{Eq:SystemControl}, the fast control $u(x,y;\xi)$ represents a direct action modifying the rate at which agents change their strategy in order to adjust the stock of renewable resource $y$ in a desired manner, for instance 
    some regulation aimed to increase or decrease such concentration by limiting the number of agents that are allowed to exploit either resource, forcing some to rapidly switch to the other consumption strategy. By contrast, the slow control $v(x,y;\xi)$ corresponds to actions directly favouring the recovery of the resource's stock, for instance, in the form of infrastructure expansions both for the storage and production of the resource, as such tasks usually require extended time intervals in order to be effectively introduced or executed.
    
    Therefore, now we present the consequent purely fast and combined fast--slow controllers for system \eqref{Eq:SystemControl}, following the results stated in Theorem \ref{Thm:GeneralizedQuadraticControl}.
    First, we consider the fast control scheme, i.e. $v(x,y;\xi)=0$, and
    
    \begin{equation}
        u(x,y;\xi) = - \frac{\varepsilon B_{c}}{\sigma k a_{c}}\left( x-x^{*} \right)(H-h)e^{-2\left(y - y^{*}\right)/\sigma \varepsilon},
        \label{Eq:FastControlMainProblem}
    \end{equation}
    with 
    \begin{equation}
        H\left( x, y,\varepsilon \right) = \frac{\sigma}{2}e^{2\left(y - y^{*}\right)/\sigma \varepsilon}\left( \frac{b_{c}}{\varepsilon}\left( y-y^{*} \right) + \frac{a_{c}}{\varepsilon}\left( x-x^{*} \right)^{2k} - \frac{\sigma b_{c}}{2} \right),
    \label{Eq:HamiltonianFunction}
    \end{equation}
    where, as explained at the beginning of this section, $a_{c}$ and $b_{c}$ are constant values obtained through the expansion of \eqref{Eq:SystemControl} near the fold point $p=(x_{p}, y_{p})\in\mathcal{C}_{0}$ of interest and are responsible of setting the parabola--like shape as well as its stability properties, while $\sigma = \text{sign}(a_{c} b_{c})$, $k\in\mathbb{N}$, $B_{c} > 0$ is the controller gain, and $(x_{p}, y_{p}) = \left( x^{*}, y^{*} \right)$, numerically obtained in our studies, are the coordinates of the fold point $p\in\mathcal{C}_{0}$ of interest in the original coordinated system of \eqref{Eq:SystemControl}, required to displace the point $p = (x_p,y_p)$ to the origin in the coordinated system of the normal form \eqref{Eq:ApproximationGeneralVectorField}. Moreover, notice that since the slow dynamics are not actuated, we are going to stabilize canards centered at $(x, y) = \left(1/r - x^{*}, 0 \right)$. Hence, it is convenient to define the coordinate transformation $\Hat{x} = x - \left(1/r -x^{*} \right)$, which brings \eqref{Eq:FastControlMainProblem} and \eqref{Eq:HamiltonianFunction} to
    \begin{equation}
        u(x,y;\xi) = -a_{c}\left( x-x^{*} \right)^{2k} + a_{c}\left( x-1/r \right)^{2k}  - \frac{\varepsilon B_{c}}{\sigma k a_{c}}\left( x-1/r \right)(H-h)e^{-2\left(y - y^{*}\right)/\sigma \varepsilon},
        \label{Eq:FastControlMainProblemFinal}
    \end{equation}
    and 
    \begin{equation}
        H\left( x, y,\varepsilon \right) = \frac{\sigma}{2}e^{2\left(y - y^{*}\right)/\sigma \varepsilon}\left( \frac{b_{c}}{\varepsilon}\left( y-y^{*} \right) + \frac{a_{c}}{\varepsilon}\left( x-1/r \right)^{2k} - \frac{\sigma b_{c}}{2} \right).
    \label{Eq:HamiltonianFunctionFinal}
    \end{equation}
    
    As a comparison to the resulting purely fast controller \eqref{Eq:FastControlMainProblemFinal} and Hamiltonian \eqref{Eq:HamiltonianFunctionFinal}, now we consider a fully actuated scenario and obtain a joint fast--slow controller for the decision--making system \eqref{Eq:SimplifiedSystem}, in the form
    \begin{align}
        \begin{split}
            u(x,y;\xi) &= -\frac{\varepsilon B_{c}}{\sigma k a_{c}}\left( x-x^{*} \right)\left( H-h \right)e^{2\left( y-y^{*} \right)/\sigma \varepsilon},\\
            v(x,y;\xi) &= -\left( 1-rx^{*} \right),
            \label{Eq:FastSlowControlMainProblem}
        \end{split}
    \end{align}
    with $H(x,y,\varepsilon)$ as \eqref{Eq:HamiltonianFunction} and the same variable interpretations as before. Notice that in the joint fast--slow scenario \eqref{Eq:FastSlowControlMainProblem}, the slow component $v(x,y;\xi)$ effectively sets the fold point $p\in\mathcal{C}_{0}$ of \eqref{Eq:SystemControl} at the origin in the coordinated system of the normal form \eqref{Eq:ApproximationGeneralVectorField}, and therefore the additional translation $\hat{x} = x-(1/r - x^{*})$ necessary for the purely fast controller \eqref{Eq:FastControlMainProblemFinal} is no longer needed. 
    
\begin{figure*}[htbp]
        \begin{subfigure}[htbp]{0.5\textwidth}
            \includegraphics[width=0.9\linewidth]{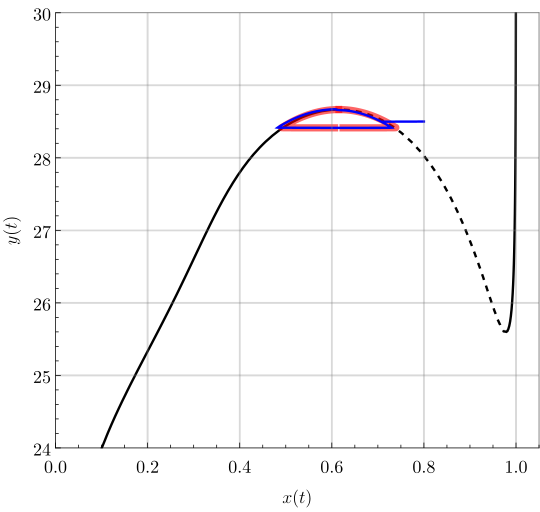}
            \end{subfigure}\hfill
            \begin{subfigure}[htbp]{0.5\textwidth}
            \includegraphics[width=0.9\linewidth]{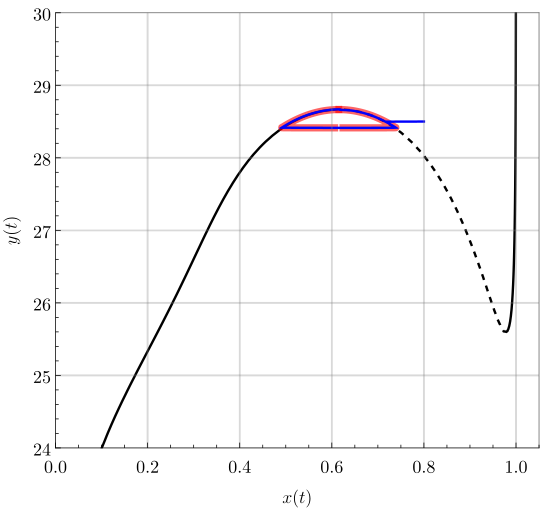}
            \end{subfigure}
            \begin{subfigure}[htbp]{0.5\textwidth}
            \includegraphics[width=0.9\linewidth]{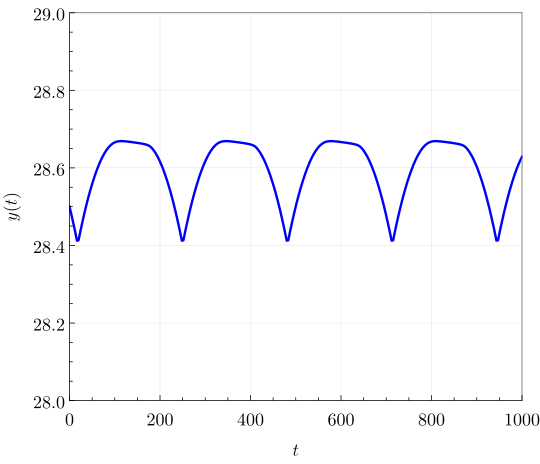}
            \end{subfigure}\hfill
            \begin{subfigure}[htbp]{0.5\textwidth}
            \includegraphics[width=0.9\linewidth]{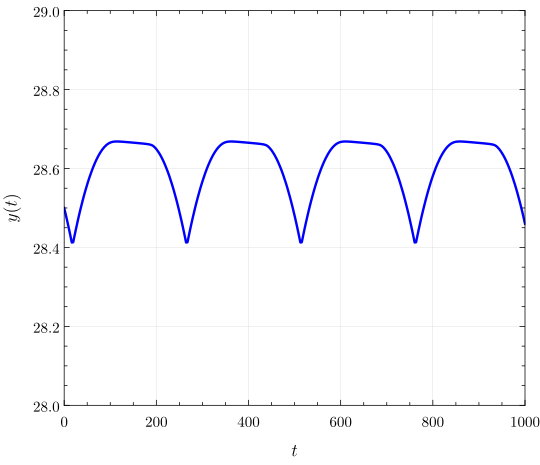}
            \end{subfigure}
            \begin{subfigure}[htbp]{0.5\textwidth}
            \includegraphics[width=0.9\linewidth]{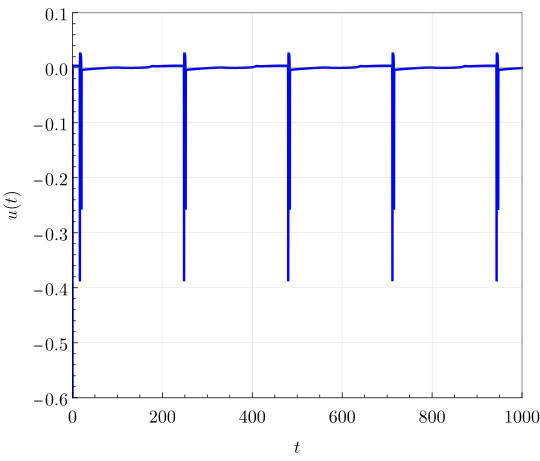}
            \end{subfigure}\hfill
            \begin{subfigure}[htbp]{0.5\textwidth}
            \includegraphics[width=0.9\linewidth]{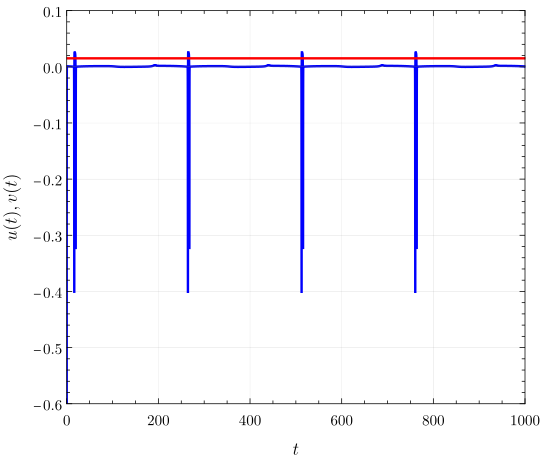}
            \end{subfigure}
            \caption{
            Purely fast \eqref{Eq:FastControlMainProblem} (left) vs.\ combined fast--slow \eqref{Eq:FastSlowControlMainProblem} (right) control in system \eqref{Eq:SimplifiedSystem}. Rows: phase portraits (top), resource stock $y(t)$ (middle), controller output (bottom). In the phase portraits, solid/dashed black are the attracting/repelling branches, with target level set $\gamma_{h}$ in red and the trajectory in blue; for the fast--slow scheme the slow input is the constant $v(x,y;\xi)=-(1-rx^{*})$ (red). Both schemes drive the trajectory onto $\gamma_{h}$ even from initial conditions near the repelling branch. Parameters: $\alpha=2.0$, $\beta=0.75$, $\gamma=0.5$, $c=2.5$, $d=1.18$, $b=30.0$, $r=1.65$, $\varepsilon=0.01$, $a_{c}=1.64218$, $b_{c}=0.100924$, $c_{c}=3.0$, $k=1$, $t=1000$, $(x(0),y(0))=(0.8,28.0)$, $(x^{*},y^{*})=(0.6163,28.665)$; gain $B_{c}=1500.0$ (left), $1000.0$ (right).}
            
            \label{Fig:CompetitiveSystemControl}
        \end{figure*}

\begin{figure*}[htbp]
        \begin{subfigure}[htbp]{0.5\textwidth}
            \includegraphics[width=0.88\linewidth]{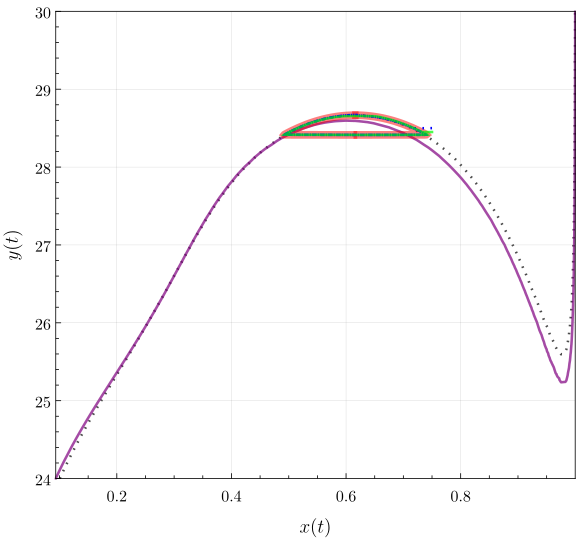}
            \end{subfigure}\hfill
            \begin{subfigure}[htbp]{0.5\textwidth}
            \includegraphics[width=0.93\linewidth]{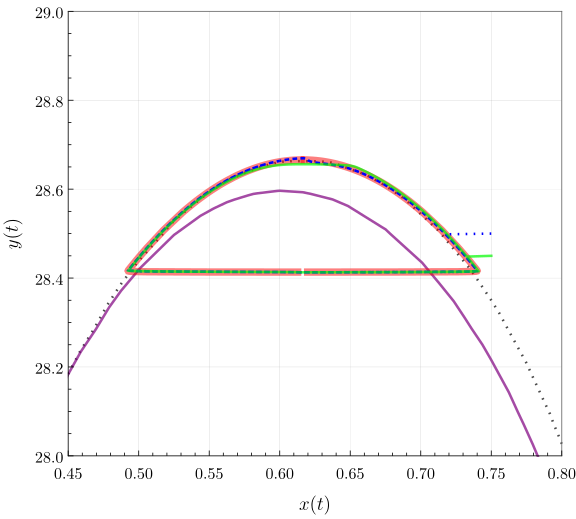}
            \end{subfigure}
            \vskip \baselineskip
            \begin{subfigure}[htbp]{0.5\textwidth}
            \includegraphics[width=0.93\linewidth]{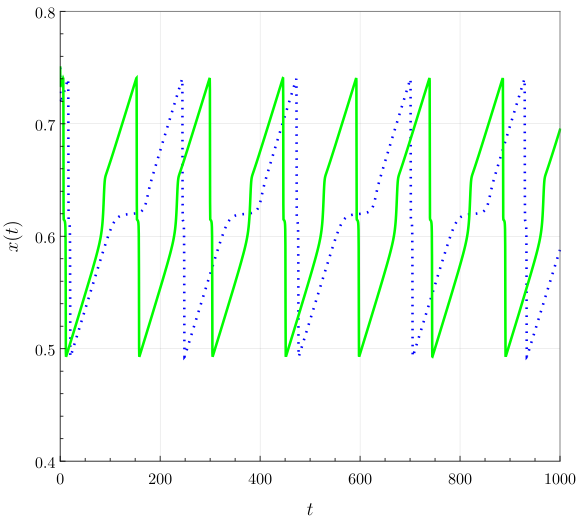}
            \end{subfigure}\hfill
            \begin{subfigure}[htbp]{0.5\textwidth}
            \includegraphics[width=0.95\linewidth]{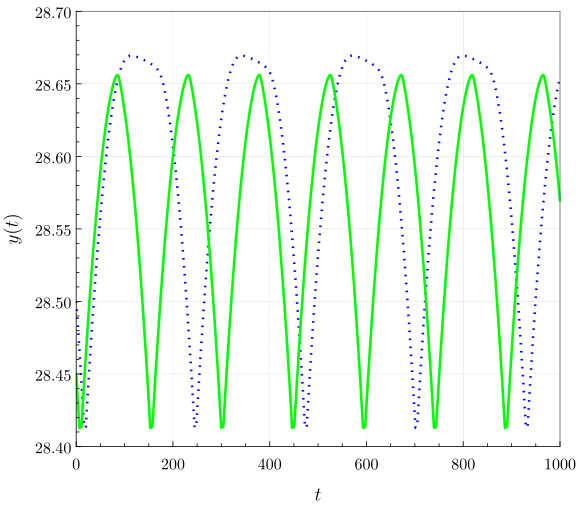}
            \end{subfigure}
            \caption{
            Stabilization of a canard cycle near $F_1$ in the decision–making system \eqref{Eq:SystemControl} under the fast–slow control \eqref{Eq:FastSlowControlMainProblem}, with (purple/green, solid) and without (black/blue, dotted) parametric perturbations; target level set $\gamma_h$ in red, $H$ as in \eqref{Eq:HamiltonianFunction}. Expansion at $F_1$: unperturbed $(a_c,b_c)=(1.6423,0.100927)$, perturbed $(a_p,b_p)=(1.75871,0.108357)$ at $(x^*,y^*)=(0.6025,28.596378)$, giving $\delta a=0.116412$, $\delta b=0.00743007$, which satisfy \eqref{Eq:ConditionPerturbationStability}. The control uses only unperturbed data yet stabilizes both cases; the lower row ($x(t)$, $y(t)$) shows only a phase shift between them.}
            \label{Fig:PerturbationDecisionMaking}
        \end{figure*}

    In Figure \ref{Fig:CompetitiveSystemControl}, we compare the results obtained by implementing the fast control \eqref{Eq:FastControlMainProblem} (left), and the fast--slow scheme \eqref{Eq:FastSlowControlMainProblem} (right) in system \eqref{Eq:SimplifiedSystem} for the stabilization of canards in a vicinity of the leftmost fold point $F_1$ (see Figure \ref{Fig:CompetitiveSystemOpenLoop}). This particular scenario represents a situation in which the controlling entity aims to increase the slowly renewable resource $y$ stock by limiting the share of agents consuming it. Observe that both controllers effectively cause trajectories traveling near a neighborhood of the fold point to follow the targeted level set, namely $\gamma_{h}$, in red, thus moving sufficiently close to unstable branches of the critical manifold and producing sustained canard orbits even when the initial conditions are set near the unstable branch of the critical manifold, showing the effectiveness of our control schemes. Particularly, notice the presence of a shift between the targeted orbit and the actual response in the purely fast controller due to the translation $\Hat{x} = x - (1/r - x^{*})$. Moreover, observe that in both cases the controllers are active only during the direction change in the limited resource stock $y(t)$, demonstrating a remarkable energetic efficiency. Now, it is clear why stabilizing the consumption ratio in the fold point $F_1$ is desirable, however, in order to precisely stabilize the point $F_1$, the controlling entity would need to accurately know the value of each parameter of the system, something that is rarely, if not never seen, in real--world applications. Hence, the controlling authority, for instance a government, aims to robustly stabilize a neighboring trajectory that stays near the fold point $F_1$, hindering any abrupt transition possible even in the presence of parametric uncertainties. To achieve the aforementioned, we extend our analysis for the stabilization of canard cycles in the perturbed normal form \eqref{Eq:PerturbedNormalForm}, detailed in Section \ref{SubSec:PerturbationNormalForm}, and the results are depicted in Figure \ref{Fig:PerturbationDecisionMaking}, where a canard cycle in a vicinity of the leftmost fold point $F_{1}$ is stabilized even in the presence of modelling perturbations by implementing the fast--slow control \eqref{Eq:FastSlowControlMainProblem} in \eqref{Eq:SystemControl}. In particular, for this example we consider perturbations $\delta_{\beta}=0.05$ and $\delta_{\gamma}=-0.001$, related to parameters $\beta$ and $\gamma$, respectively. In the phase portraits shown in the upper row, we represent the critical manifold and response for the unperturbed system \eqref{Eq:SystemControl} in dotted black and blue, while the critical manifold and response for the same system considering the perturbations $\delta_{\beta}$ and $\delta_{\gamma}$ are presented in solid purple and green, correspondingly. Additionally, the target level set $\gamma = \{ (x,y)\in \mathbb{R}^{2}: H=h \}$ is presented in red, where $H(x,y,\varepsilon)$ is the Hamiltonian \eqref{Eq:HamiltonianFunction} of the unperturbed model \eqref{Eq:SystemControl}. Observe that the canard orbit is effectively stabilized for both the unperturbed and perturbed scenarios. Specifically, the parameters of the parabola--like shape obtained through the expansion of \eqref{Eq:SystemControl} around the leftmost fold point $F_{1}$ for the unperturbed case are $a_{c}=1.6423$ and $b_{c}=0.100927$, while for the perturbed scenario are $a_{c}=1.75871$ and $b_{c}=0.108357$, with the fold point $F_{1}$ of the perturbed case centered at $(x^{*}, y^{*}) = (0.6025, 28.596378)$, once again numerically identified. Therefore, the resultant perturbations, given by the parameter difference, are $\delta_{a}=0.116412$ and $\delta_{b}=0.00743007$, satisfying condition \eqref{Eq:ConditionPerturbationStability}. Additionally, in the lower row the same phase--shift observed for the perturbed normal form \eqref{Eq:PerturbedNormalForm} is appreciated in the decision--making system \eqref{Eq:SystemControl} under perturbations. Nevertheless, the qualitative behaviour of \eqref{Eq:SystemControl} with and without perturbations are similar. Finally, observe that the origin of the target level set $\gamma_{h}$ is centered at the leftmost fold point $F_{1}$ of the unperturbed scenario as the controller is designed only with information of the unperturbed case, granting the controlling entity with the capability to stabilize orbits in a vicinity of the desired fold point even with a bounded level of parametric uncertainty.

In the following section, we further explore the capabilities of our fast--slow controllers \eqref{Eq:FastControlMainProblemFinal} and \eqref{Eq:FastSlowControlMainProblem} to stabilize canard cycles in a vicinity of the rightmost fold point $F_2$ of the decision--making model \eqref{Eq:SimplifiedSystem} and discuss on its utility in our system.

    \subsection{Control of the optimal consumption strategy}
    \label{SubSec:OptimalStrategy}
\begin{figure*}[htbp]
        \begin{subfigure}[htbp]{0.5\textwidth}
            \includegraphics[width=0.88\linewidth]{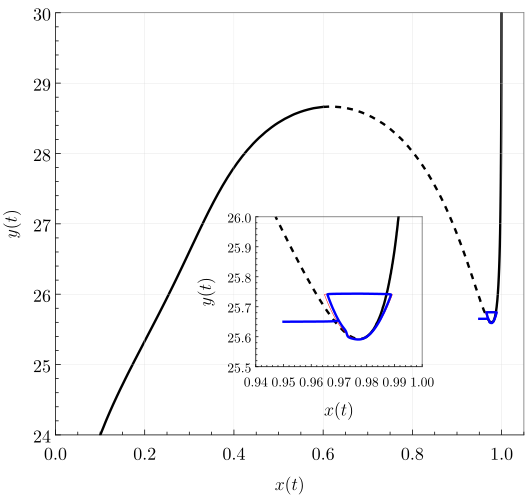}
            \caption{}
            \end{subfigure}\hfill
            \begin{subfigure}[htbp]{0.5\textwidth}
            \includegraphics[width=0.97\linewidth]{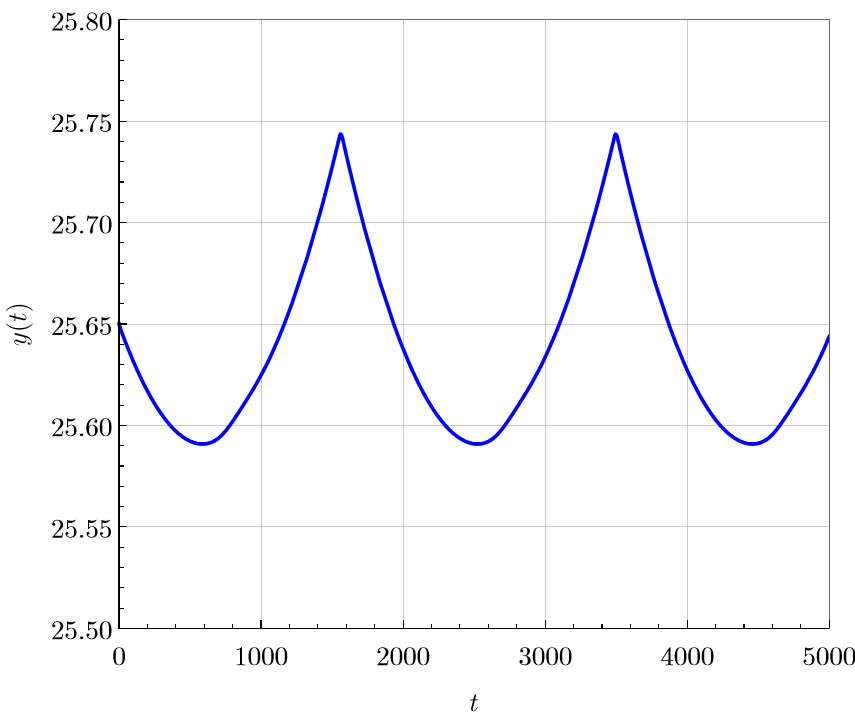}
            \caption{}
            \end{subfigure}
            \vskip \baselineskip
            \begin{subfigure}[htbp]{0.5\textwidth}
            \includegraphics[width=0.93\linewidth]{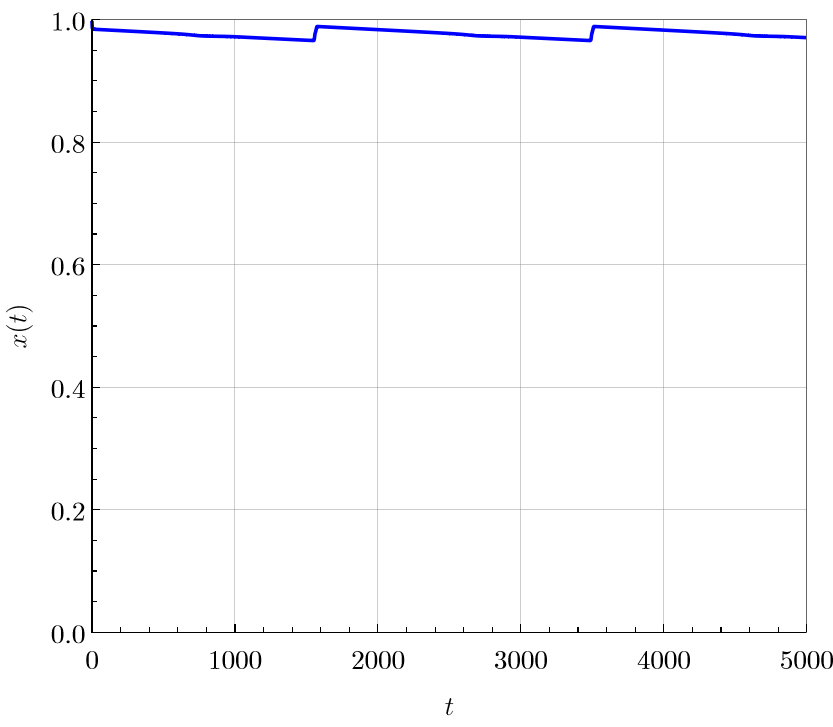}
            \caption{}
            \end{subfigure}\hfill
            \begin{subfigure}[htbp]{0.5\textwidth}
            \includegraphics[width=0.97\linewidth]{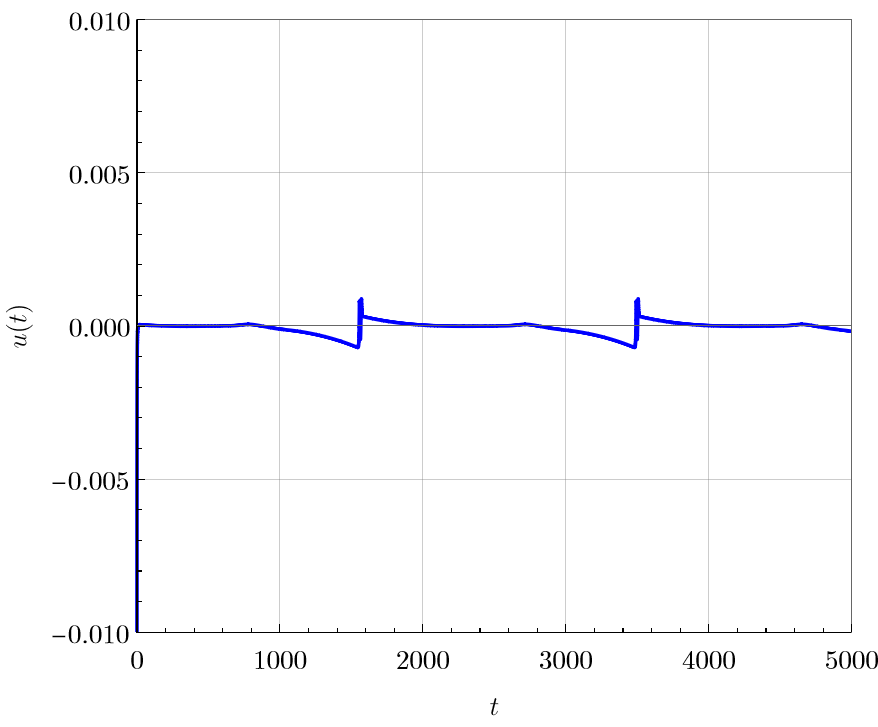}
            \caption{}
            \end{subfigure}
            \caption{Stabilization of canard trajectories around the optimal consumption fold point $F_2$. \textbf{a)} $x-y$ plane, with the stable and unstable branches of the critical manifold in solid and dashed black, respectively, while the target level set $\gamma_{h}$ and the resulting trajectory are shown in red and blue, correspondingly. \textbf{b)} Variations in the slowly renewable resource's stock $y(t)$, \textbf{c)} agents consuming the limited resource, and \textbf{d)} fast--slow controller response for the stabilization of the rightmost fold $F_2$.}
            \label{Fig:RightmostFoldF2}
        \end{figure*}

\begin{figure*}[htbp]\centering
            \begin{subfigure}[htbp]{0.5\textwidth}
                \includegraphics[height=0.8\textwidth]{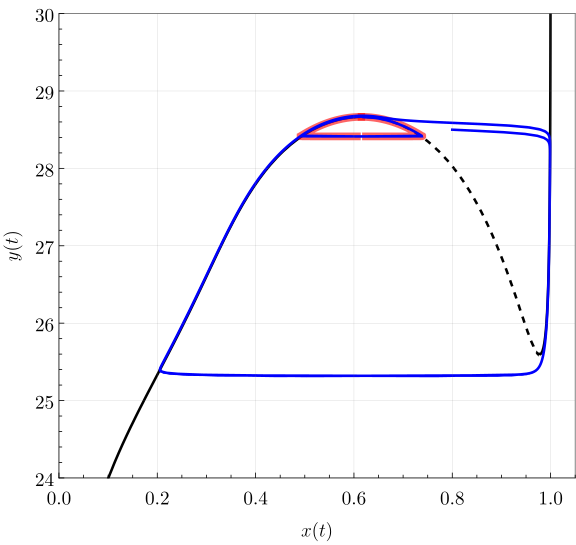}
                \end{subfigure}\hfill
                \begin{subfigure}[htbp]{0.5\textwidth}
                \includegraphics[height=0.8\textwidth]{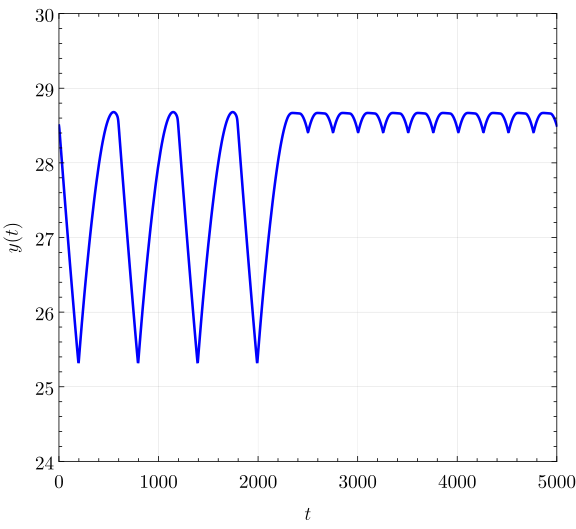}
                \end{subfigure}
                \begin{subfigure}[htbp]{0.5\textwidth}
                \includegraphics[height=0.8\textwidth]{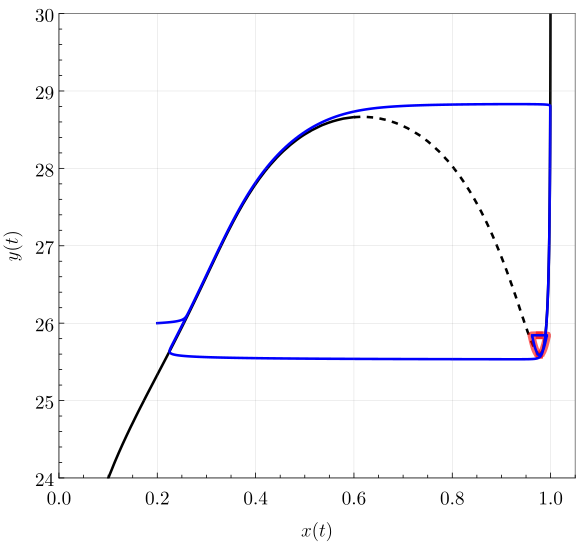}
                \end{subfigure}\hfill
                \begin{subfigure}[htbp]{0.5\textwidth}
                \includegraphics[height=0.8\textwidth]{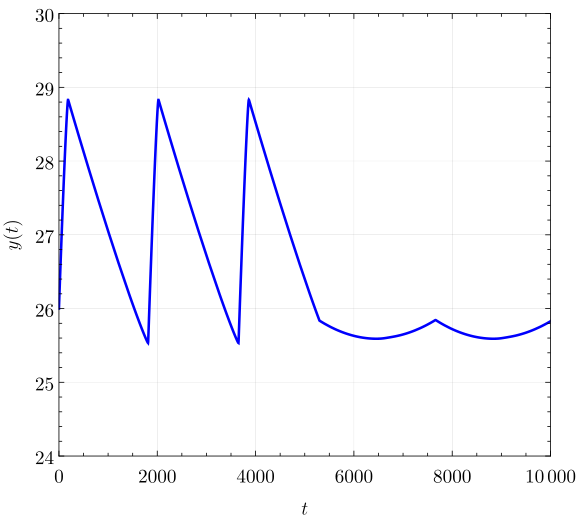}
                \end{subfigure}
                \begin{subfigure}[htbp]{0.5\textwidth}
                \includegraphics[height=0.8\textwidth]{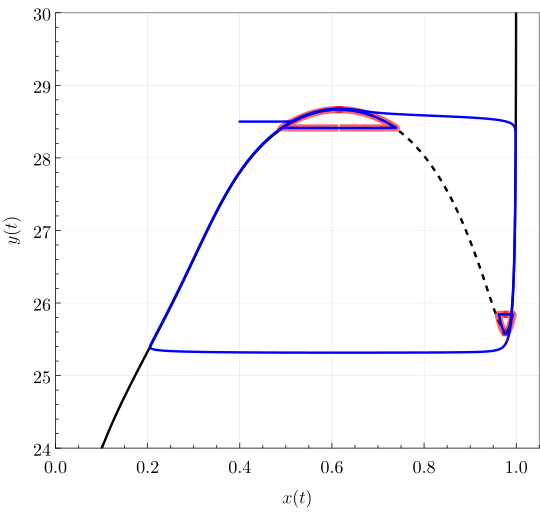}
                \end{subfigure}\hfill
                \begin{subfigure}[htbp]{0.5\textwidth}
                \includegraphics[height=0.8\textwidth]{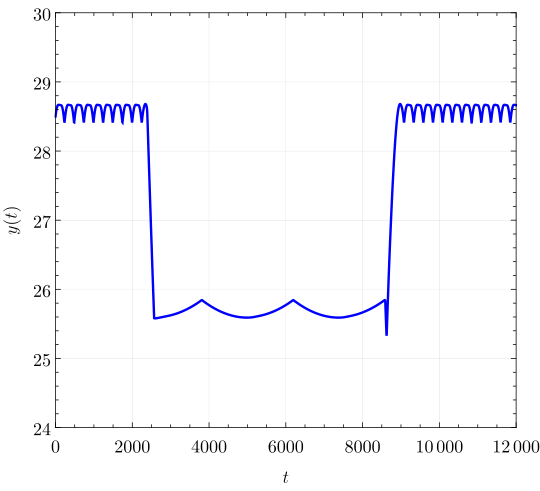}
                \end{subfigure}
                \caption{
                Activation of the fast--slow control \eqref{Eq:FastSlowControlMainProblem} on the decision--making model \eqref{Eq:SimplifiedSystem}, stabilizing a canard orbit near $F_1$ (upper), $F_2$ (middle), and sequentially both folds (lower). Left column: $x$--$y$ planes (solid/dashed black: attracting/repelling branches; target sets $\gamma_{h}$ in red, response in blue). Right column: the resource $y$ over time. Common parameters: $\alpha=2.0$, $\beta=0.75$, $\gamma=0.5$, $b=30$, $c=2.5$, $d=1.18$, $\varepsilon=0.001$. $F_1$ (upper): $r=1.62$, $(x(0),y(0))=(0.8,28.5)$, $(x^{*},y^{*})=(0.6163,28.665)$, $t_{on}=2300$. $F_2$ (middle): $r=1.08$, $(x(0),y(0))=(0.2,26.0)$, $(x^{*},y^{*})=(0.977,25.5928)$, $t_{on}=5300$. Sequential (lower): same parameters, switch times $t_{off_{F_1}}=2300$, $t_{on_{F_2}}=2570$, $t_{off_{F_2}}=8590$, $t_{on_{F_1}}=9000$.}
                \label{Fig:ActivationControlOnOff}
            \end{figure*}

        Lastly, we extend the implementation of our control schemes in the decision--making model \eqref{Eq:SimplifiedSystem} to the stabilization of a canard trajectory in the vicinity of the rightmost fold point, namely $F_{2}$ (see Figure \ref{Fig:CompetitiveSystemOpenLoop}). The selection of $F_2$ is related to an \emph{optimal exploitation strategy} as this point represents the limit at which the largest group of agents benefit from consuming the slowly renewable resource $y$, while  also being the point at which the minimum stock of resource $y$ is reserved, advantageous from an economic perspective since it considerably reduces the storage costs for the controlling entity. On top of that, the 
        controlling authority is entitled once more with the capability to produce abrupt transitions by exploiting the system's criticality, for instance in a scenario in which now a greater concentration of resource in needed. As an example of the implementation of our fast--slow controllers, in Figure \ref{Fig:RightmostFoldF2} the effective stabilization of canard cycles in a vicinity of the rightmost fold point $F_2$ is shown along with the resulting resource, agents and control responses when using the fast--slow control scheme \eqref{Eq:FastSlowControlMainProblem}. Observe that the variation of agents consuming the renewable resource is rather small, while the element $y$ also presents periodic oscillations with small amplitude. In addition, notice that the controller only activates periodically and remains bounded between two considerably small values, which translates in minor actions done by the controlling entity in order to stay around the desired trajectory and optimally consume the renewable resource $y$. Finally, in Figure \ref{Fig:ActivationControlOnOff} we show a scenario in which the controlling authority effectively regulates at will the dynamical behaviour of \eqref{Eq:SimplifiedSystem} by activating our fast--slow control scheme \eqref{Eq:FastSlowControlMainProblem}. On the upper row of Figure \ref{Fig:ActivationControlOnOff}, we show the stabilization of an orbit near the leftmost fold $F_1$, which corresponds to the scenario when the share of agents following each consumption strategy is mostly balanced, leading to the largest amount of renewable resource $y$ in stock, without considering the right asymptote. From an authority perspective, stabilizing a canard cycle in a vicinity of $F_1$ represents an advantageous position in order to increase the amount of limited resource, to then allow its consumption by a larger share of agents after a specific desired time, corresponding to the free dynamics of \eqref{Eq:SimplifiedSystem}. On the other hand, in the middle row of Figure \ref{Fig:ActivationControlOnOff} we present the effect of activating our controller after a given time for the stabilization of a canard cycle in a vicinity of the fold point $F_2$, which corresponds to the aforementioned optimal consumption strategy. Similarly, in the lower row of Figure \ref{Fig:ActivationControlOnOff} we show a sequential activation of our fast--slow controllers in order to generate a pattern that oscillates in a vicinity of the leftmost fold $F_{1}$, then deactivate the controller and let the system travel freely near the critical manifold, to later activate the controller once more when the response is in a vicinity of the rightmost fold $F_{2}$ to stabilize the canard orbit around a neighboring target level set, deactivate again the controller and finally activate it again in a vicinity of the leftmost fold point $F_{1}$, showing the capability of the controlling entity to exploit the systems' criticality to reach the desired state by sequentially activating and deactivating the controller \eqref{Eq:FastSlowControlMainProblem}. Notice that the sharp transition in the stock of renewable resource $y$, occurring while changing the target orbit from $F_{2}$ to $F_{1}$, is due to the own flow of the open--loop problem \eqref{Eq:SimplifiedSystem}. Finally, observe that the oscillation frequency for the canard orbits stabilized near the fold point $F_1$ is comparably larger with respect to the one of the canard orbits stabilized in a neighborhood of the fold point $F_2$. This fact is due to the form of the critical manifold $\mathcal{C}_0$, as the resulting trajectory mainly travels along the slow direction near the fold point $F_2$, while in the case for a canard cycle close the fold $F_1$ the systems' response constantly alternates between the fast and slow directions.

\section{Conclusion and Discussion}\label{Sec:Discussion}
Decision--making represents a fundamental component in several real world phenomena as it describes a common problem in which agents need to find the maximum benefit according to different external and internal factors. In this work, we have studied a fast--slow dynamical system describing the decision--making process that two groups of clearly identified agents undergo when selecting between two harvesting strategies, one representing the exploit of an unlimited but highly costly common source, and the other describing the consumption of a comparably cheaper but limited and slowly renewable resource. With this in mind and by an extensive numerical study, we determine parameter combinations such that the associated critical manifold of the fast--slow dynamical system analyzed presents either zero, two or four fold points, allowing for the generation of canard cycles due to the system's properties. Such special trajectories are well--known for following unstable regions of the slow--flow for considerable amounts of time, convenient in our context for the controlling entity to exploit the systems' criticality in order to manipulate the renewable resource dynamics according to a desired outcome and enabling the possibility to produce abrupt strategic transitions in the consumption schemes that the agents follow. Therefore, we have been able to stabilize canard orbits in our system by implementing fast--slow controllers developed through the analysis of a canonical form of a not necessarily generic fold point, and shown the effectiveness of our approach for two particular scenarios when considering a parameter combination with two fold points. The first one being a situation in which the controlling authority needs to guarantee a high concentration of renewable resource for a specific period, task achieved by implementing policies that restrict the number of agents consuming such source and enabling the controlling entity to cause an abrupt transition to a scenario in which most of the population consume the limited resource with a minor control effort. On the other hand, the second implementation describes a circumstance in which, lead by the controlling entity, the population follows an optimal consumption strategy that allows the vast majority of the entire community to consume the renewable resource while reserving a local minimum of it, which represents an economic saving considering the stock's storage. Moreover, we demonstrate that our fast--slow controllers present a certain degree of robustness to parametric perturbations, which grants the controlling entity with the capability to stabilize canard orbits even in the presence of modelling uncertainties.

     However, we want to emphasize that our results are able to stabilize canard orbits in any fast--slow system that presents planar folded canard orbits, regardless on the degeneracy degree such bifurcation points may present, allowing for a wide variety of real--world processes to benefit from our controllers. For instance, consider competition, which represents a key element in the study of diverse phenomena, ranging from biological systems \cite{Ref:Abrams2022}, epidemics \cite{Ref:Ye2022}, ecology \cite{Ref:Zeppini2020, Ref:Chen2021}, economics \cite{Ref:Asano2021, Ref:Yao2022}, and social models \cite{Ref:Bravetti2018}, describing the interaction between populations of a broad nature, behaviour that considerably resembles the decision--making process that we have studied.
     
Our analysis provides several future research directions:
\begin{enumerate}
  \item[(i)] \emph{Coupled decision--making units.} When several populations of the form
  \eqref{eq:main2} are coupled, through a shared resource stock or social interaction between
  their shares $x_i$, the natural question is the synchronisation of their canard--mediated
  transitions, and whether the controllers of Section~\ref{Sec:Results} can enforce or break
  it. The target canards need not be of the same type: coordinating units at the same fold
  differs from a crossed configuration ($F_1$--$F_2$), where the coupled canards have
  different amplitude and frequency.
  \item[(ii)] \emph{Multiple limited resources.} The slow space of \eqref{eq:main2} is
  one--dimensional because a single limited resource is tracked; strategies enter only through
  the share $x$. With several limited, slowly renewable resources, each contributing its own
  stock, the slow space becomes higher--dimensional, providing the possibility of
  folded--node canards and mixed--mode oscillations \cite{Ref:Szmolyan2001,Ref:Desroches2016}. Extending
  Theorem~\ref{Thm:GeneralizedQuadraticControl}, which depends on the scalar Hamiltonian
  \eqref{Eq:HamiltonianFunctionFinal}, to this setting, where no such Hamiltonian is available, is open.
  \item[(iii)] \emph{Stochastic canard control.} Bounded--rational exploration makes noise
  intrinsic to the decision dynamics. As the canard window has exponentially small width
  $\mathcal{O}(e^{-K/\varepsilon})$, weak noise can dominate it and trigger unintended
  transitions; whether the feedback \eqref{Eq:FastControlEntry} survives such forcing or must
  be redesigned is open.
  \item[(iv)] \emph{Mechanism of the change of criticality.} Figure~\ref{fig:criticalitymaps}
  shows the criticality invariant $\sigma$ of the Hopf bifurcation changing sign, separating a supercritical from a subcritical onset of oscillations. In turn, this determines whether a policy past the bifurcation
  is reversible. Identifying the mechanism, a Bautin point or a fold collision, and locating
  the codimension--two points analytically is open, since $f$ admits no closed--form Lyapunov
  coefficient.
  \item[(v)] \emph{Robust global control and fold identification.} Our controllers are local,
  assuming known fold coordinates $(x^{\ast},y^{\ast})$ and an accurate parabola--like
  approximation. A switching scheme incorporating a manifold--following controller once the
  local error grows, and the online estimation of $(x^{\ast},y^{\ast})$ for adaptive
  stabilisation, are natural extensions.
\end{enumerate}

\appendix
\renewcommand{\thetheorem}{\Alph{section}.\arabic{theorem}}
\renewcommand{\thedefinition}{\Alph{section}.\arabic{definition}}
\renewcommand{\theremark}{\Alph{section}.\arabic{remark}}
\section{Background on slow-fast systems}\label{Sec:Preliminaries}
    In this section, we present the theoretical principles that represent the basis of our results. We begin by briefly discussing the theory of fast--slow dynamics and geometric singular perturbation theory, specifically Fenichel's theorem \cite{Ref:Fenichel1979, Ref:Kuehn2015}. Thereafter, we describe planar folded canard trajectories, emphasizing their normal form \cite{Ref:Krupa2001, Ref:JardonKojakhmetov2022}.
\subsection{Fast--slow systems}
    \label{Sec:FastSlowSystems}
        Consider the class of systems given by
        \begin{align}
            \begin{split}
            \varepsilon \dot{x} &= f(x,y,\varepsilon),\\
            \dot{y} &= g(x,y,\varepsilon),
            \label{Eq:FastSlowSystem_Slow}
            \end{split}
        \end{align}
        typically known as a \emph{fast--slow vector field}, composed by a set of singularly perturbed ordinary differential equations, where the over--dot denotes the derivative with respect to the slow--time $\tau$, and with $x=x(\tau)\in\mathbb{R}^{m}$, $y=y(\tau)\in\mathbb{R}^{n}$ the fast and slow variables, respectively. Moreover, the parameter $0<\varepsilon\ll1$ describes the timescale separation and the mappings $f:\mathbb{R}^{m}\times \mathbb{R}^{n} \times \mathbb{R} \xrightarrow{} \mathbb{R}^{m}$ and $g: \mathbb{R}^{m} \times \mathbb{R}^{n} \times \mathbb{R}\xrightarrow{} \mathbb{R}^{n}$ are assumed to be sufficiently smooth. Additionally, by defining the fast time variable $t\coloneqq \tau/\varepsilon$, it is possible to express the equivalent fast--time form of \eqref{Eq:FastSlowSystem_Slow} as
        \begin{align}
            \begin{split}
                x' &= f(x,y,\varepsilon),\\
                y' &= \varepsilon g(x,y,\varepsilon),
                \label{Eq:FastSlowSystem_Fast}
            \end{split}
        \end{align}
        with the prime denoting the derivative with respect to $t$. Now, one important mathematical theory usually employed in the analysis of systems \eqref{Eq:FastSlowSystem_Slow}--\eqref{Eq:FastSlowSystem_Fast} is geometric singular perturbation theory (GSPT), with an overall idea of analyzing both problems separately and looking for invariant objects that can persist certain small perturbations \cite{Ref:Kuehn2015}. Specifically, by considering the singular limit $\varepsilon = 0$ in \eqref{Eq:FastSlowSystem_Slow} and 
        \eqref{Eq:FastSlowSystem_Fast}, two different problems arise
        \begin{equation}
            \begin{aligned}
            & 0=f(x,y,0), \qquad \qquad \qquad && x'=f(x,y,0),\\
            & \dot{y} = g(x,y,0), \qquad \qquad \qquad && y'=0,
            \label{Eq:SingularLimits}
            \end{aligned}
        \end{equation}
        known, correspondingly, as the reduced slow subsystem, being a constrained differential equation \cite{Ref:Takens1976} or a differential algebraic equation \cite{Ref:Kunkel2006}, and the layer problem. Naturally, the aforementioned problems are no longer equivalent but they are intrinsically related through the geometric object known as the \emph{critical manifold}.
        \begin{definition}[Critical manifold]
            The critical manifold of a fast--slow problem is
            \begin{equation}
                \mathcal{C}_{0} = \left\{ (x,y)\in\mathbb{R}^{m}\times \mathbb{R}^{n} : f(x,y,0) = 0 \right\}.
                \label{Eq:CriticalManifoldDefinition}
            \end{equation}
        \end{definition}

        \begin{remark}
            Observe that the critical manifold \eqref{Eq:CriticalManifoldDefinition} corresponds both to the phase--space of the reduced problem and the set of equilibrium of the layer problem.
        \end{remark}

        Moreover, one highly relevant property of fast--slow systems to analyze the different types of possible bifurcations is the \emph{normal hyperbolicity}, defined as follows.

        \begin{definition}[Normal hyperbolicity]\label{def:nh}
            A subset $\mathcal{S}_{0} \subset \mathcal{C}_{0}$ is normally hyperbolic if the matrix $(\text{D}_{x}f)(p,0)$ has only eigenvalues with non--zero real part for every point $p\in \mathcal{S}_{0}$. Furthermore, given a normally hyperbolic subset $\mathcal{S}_{0}$, it can be attracting (repelling) if every eigenvalue of the Jacobian $(\text{D}_{x}f)(p,0)$ has negative (positive) real part for every $p\in\mathcal{S}_{0}$. Finally, if $\mathcal{S}_{0}$ is normally hyperbolic but neither attracting nor repelling, then it is of saddle type.
        \end{definition}

A useful consequence of normal hyperbolicity is that the dynamics of the constraint equation in \eqref{Eq:SingularLimits} can be reduced to the dynamics of an ODE. Namely, due to normal hyperbolicity, the critical manifold $\mathcal{C}_{0}$ can be given locally near a regular point $p\in\mathcal{C}_{0}$ as the graph of a function $f(h(y),y,0)=0$, with $h:\mathbb{R}^{n}\mapsto\mathbb{R}^{m}$. Thus, it is possible to reduce the slow subsystem
\begin{equation}
    \begin{aligned}
        0&=f(x,y,0),\\
        \dot{y}&=g(x,y,0),
    \end{aligned}
\end{equation}
to the simpler form
\begin{equation}
    \dot{y}=g(h(y),y,0).
\end{equation}
        Hence, a subset $\mathcal{S}_{0}$ is non--normally hyperbolic if the Jacobian $(\text{D}_{x}f)(p,0)$ has at least one eigenvalue with zero real part. According to the normally hyperbolic characteristic of the problem at hand, different analysis techniques can be used. For instance, non--normally hyperbolic points, which are related to dynamic features such as relaxation oscillations, that is, limit cycles of a singularly perturbed dynamical system \cite{Ref:Grasman2011}, appearing in chemistry \cite{Ref:Rossler1978}, geophysics \cite{Ref:Vasconcelos1996}, biology \cite{Ref:Deng2004, Ref:Shen2016, Ref:Zhao2022}, or economics \cite{Ref:Hamburger1934, Ref:Varian1979}, and canard trajectories \cite{Ref:Desroches2016}, can be studied with the blow--up method \cite{Ref:Dumortier1991, Ref:Alvarez2011, Ref:Kuehn2015, Ref:Jardon2021}. On the other hand, for normally hyperbolic cases a widely employed technique comes from the Fenichel's theorem, stated as follows.

        \begin{theorem}[Fenichel's theorem]
            Suppose $\mathcal{S}=\mathcal{S}_{0}$ is a compact normally hyperbolic submanifold of the critical manifold $\mathcal{C}_{0}$ of \eqref{Eq:FastSlowSystem_Slow}, and that $f$, $g\in C^{r}$, $r<\infty$. Then, for  $0<\varepsilon\ll1$ sufficiently small, the following hold:
            \begin{enumerate}
                \item There exists a locally invariant manifold $\mathcal{S}_{\varepsilon}$ diffeomorphic to $\mathcal{S}_{0}$. In this scope, local invariance means that trajectories can enter or escape $\mathcal{S}_{\varepsilon}$ only through its boundaries.
                \item The Hausdorff distance between the submanifolds $\mathcal{S}_{\varepsilon}$ and $\mathcal{S}_{0}$ is of order $\varepsilon$.
                \item The flow on $\mathcal{S}_{\varepsilon}$ converges to the flow generated by the reduced problem \eqref{Eq:FastSlowSystem_Slow} in the singular limit $\varepsilon=0$.
                \item $\mathcal{S}_{\varepsilon}$ is $C^{r}$ smooth, with $r<\infty$.
                \item $\mathcal{S}_{\varepsilon}$ shares the same normally hyperbolic and stability properties as the critical manifold $\mathcal{S}_{0}$.
                \item $\mathcal{S}_{\varepsilon}$ is rarely unique, as all submanifolds satisfying properties $1-5$, in a region with fixed distance from $\partial \mathcal{S}_{\varepsilon}$, are exponentially close to each other with a Hausdorff distance of order $O(\emph{exp}(-C/\varepsilon))$, $C>0$, and $C\in O(1)$ as $\varepsilon \rightarrow 0$. Therefore, every submanifold $\mathcal{S}_{\varepsilon}$ satisfying conditions $1-5$ is called a \emph{slow manifold}.
            \end{enumerate}
            \label{Thm:Fenichel}
        \end{theorem}

        Once we have revisited the main concepts of fast--slow systems, in what follows we focus our interest in the unexpected canard trajectories, solutions that follow unstable branches of the critical manifold $\mathcal{C}_{0}$ for considerable time, as the ones shown in the right panel of Figure \ref{Fig:CompetitiveSystemOpenLoop}, and which can arise from bifurcations of the fold type. Such trajectories have been observed in various applied sciences, for instance in neuroscience, allowing a detailed description of the very fast onset of large amplitude oscillations caused by small parameter variations in neural models \cite{Ref:Izhikevich2006, Ref:Ermentrout2009, Ref:Ermentrout2010, Ref:Koskal2021}, as well as in chemistry \cite{Ref:Brons1991, Ref:Milik2001, Ref:Moehlis2002}, to mention a few.

A particularly important phenomenon near fold points of the critical manifold is the \emph{singular Hopf bifurcation}~\cite{Ref:Krupa2001,Ref:Kuehn2015}. Consider a planar fast--slow system~\eqref{Eq:FastSlowSystem_Slow} with a non-degenerate fold point $p \in C_0$ satisfying the standard conditions: $\partial_x f(p,0) = 0$, $\partial_{xx} f(p,0) \neq 0$, $\partial_y f(p,0) \neq 0$, and $g(p,0) \neq 0$. Suppose additionally that the system possesses an equilibrium point that, as a parameter (say~$\mu$) varies, crosses the fold point~$p$. For $0 < \varepsilon \ll 1$, Krupa and Szmolyan~\cite{Ref:Krupa2001} showed that this crossing produces a Hopf bifurcation of the full system at a parameter value $\mu = \mu_H(\varepsilon)$, with $\mu_H \to \mu_0$ as $\varepsilon \to 0$, where $\mu_0$ corresponds to the equilibrium lying exactly at the fold. This is called a \emph{singular} Hopf bifurcation because it degenerates in the singular limit.
 
The singular Hopf bifurcation is accompanied by a \emph{canard explosion}: as~$\mu$ increases past~$\mu_H(\varepsilon)$, the small-amplitude limit cycle born at the Hopf point grows rapidly through a sequence of canard cycles until it reaches the size of a full relaxation oscillation \cite{Ref:Kuehn2015}. This entire transition occurs within a parameter interval of exponentially small width $O(e^{-c/\varepsilon})$, $c > 0$. The canard cycles within this explosion are level sets of a Hamiltonian that arises in the blown-up analysis near the fold (see~\ref{Sec:PlanarFoldedCanards}), a fact that is central to the controller design in Section~\ref{Sec:Results}.

The exponential sensitivity of canard cycles to parameter variations is the fundamental reason why these trajectories, despite being dynamically significant, are difficult to observe in uncontrolled systems subject to even minor disturbances. This motivates the development of feedback controllers that stabilise canard cycles as robust attractors, which is the main objective of the present work.

    \subsection{Planar folded canards}
    \label{Sec:PlanarFoldedCanards}
\begin{figure}[h!]
            \centering
            \includegraphics[width=0.5\linewidth]{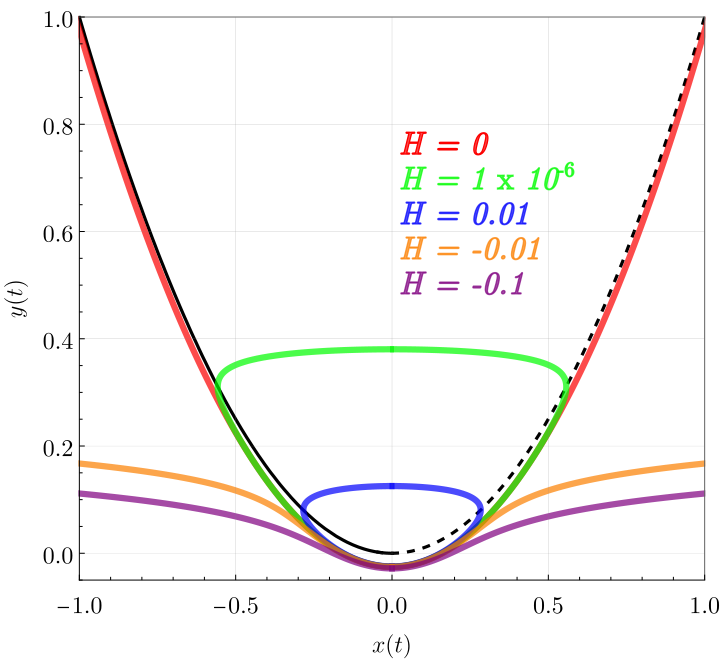}
            \caption{Canard trajectories of system \eqref{Eq:NormalFormNonHigherOrderTerms} obtained as level sets of the Hamiltonian \eqref{Eq:HamiltonianNormalForm}, where the attracting and repelling regions of the critical manifold $\mathcal{C}_{0}$ are depicted in solid and dashed black, respectively. Observe that the orange ($H=-0.01$) and purple ($H=-0.1$) solutions describe canard trajectories, while the blue ($H=0.01$) and green ($H=1\times 10^{-6}$) are canard cycles. Finally, the red orbit $(H=0)$ represents the maximal canard, which connects the attracting and repelling slow manifolds of \eqref{Eq:NormalFormNonHigherOrderTerms}.}
            \label{Fig:CanardsExample}
        \end{figure}

         Let us start by recalling the canonical form of a canard point \cite{Ref:Krupa2001}, given as
        \begin{align}
            \begin{split}
                \dot{x} &= -yh_{1}\left( x, y, \varepsilon, \alpha \right) + x^{2}h_{2}\left( x, y, \varepsilon, 
                \alpha\right) + \varepsilon h_{3}\left( x, y, \varepsilon, \alpha \right),\\
                \dot{y} &= \varepsilon \left( xh_{4}\left( x, y, \varepsilon, \alpha \right) - \alpha h_{5}\left( x, y, \varepsilon, \alpha \right) + y h_{6}\left( x, y, \varepsilon, \alpha \right) \right),
            \end{split}
            \label{Eq:CanardNormalForm}
        \end{align}
        with $(x, y)\in \mathbb{R}^{2}$, $0<\varepsilon \ll 1$, and $\alpha$ is a parameter. Furthermore,
        \begin{align}
            \begin{split}
                h_{3}(x, y, \varepsilon, \alpha) &= \mathcal{O}(x, y, \varepsilon, \alpha),\\
                h_{i}(x, y, \varepsilon, \alpha) &= 1 + \mathcal{O}(x, y, \varepsilon, \alpha), \quad i = 1, 2, 4, 5,
            \end{split}
            \label{Eq:HFunctionsNormalForm}
        \end{align}
        and $h_{6}$ is smooth. Let us simplify \eqref{Eq:CanardNormalForm} along \eqref{Eq:HFunctionsNormalForm} as
        \begin{align}
            \begin{split}
                \dot{x} &= -y + x^{2} + \tilde{f}(x, y, \varepsilon, \alpha),\\
                \dot{y} &= \varepsilon(x - \alpha + \tilde{g}(x, y, \varepsilon, \alpha)),
            \end{split}
            \label{Eq:SimplifiedNormalForm}
        \end{align}
        where the functions $\tilde{f}$ and $\tilde{g}$ collect the higher order terms in \eqref{Eq:CanardNormalForm}. Thus, locally the critical manifold is denoted by
        \begin{equation}
            \mathcal{C}_{0} = \left\{ (x, y) \in \mathbb{R}^{2}: -y + x^{2} + \tilde{f}(x, y, 0, \alpha) = 0 \right\}.
            \label{Eq:CriticalManifoldNormalForm}
        \end{equation}
        In particular, considering \eqref{Eq:SimplifiedNormalForm} in the absence of higher order terms yields
        \begin{align}
            \begin{split}
                \dot{x} &= -y + x^{2},\\
                \dot{y} &= \varepsilon(x-\alpha).
            \end{split}
            \label{Eq:NormalFormNonHigherOrderTerms}
        \end{align}
        It is straightforward to check that the orbits of \eqref{Eq:NormalFormNonHigherOrderTerms}, for $\varepsilon>0$ and $\alpha = 0$, are given by the level sets of the Hamiltonian
        \begin{equation}
            H(x, y, \varepsilon) = \frac{1}{2}e^{-2y/\varepsilon}\left( \frac{y}{\varepsilon} - \frac{x^{2}}{\varepsilon} + \frac{1}{2} \right).
            \label{Eq:HamiltonianNormalForm}
        \end{equation}
        Furthermore, some orbits of \eqref{Eq:HamiltonianNormalForm} are canard trajectories, as the ones shown in Figure \ref{Fig:CanardsExample}, and in fact it is well--known that canard cycles exist for $H\in\left( 0, 1/4 \right)$ \cite{Ref:Krupa2001}. Therefore, the representation of the normal form for a folded canard is of high importance for our study as the control schemes designed are based on such description. However, it is important to highlight that one should be careful while designing such controllers as the dynamical structure of \eqref{Eq:CanardNormalForm} could be altered, rendering inaccurate behaviours. For this reason, we present the following definitions.

        \begin{definition}[$k$--jet equivalence]
        \label{Def:KJetEquivalence}
            Let $F: \mathbb{R}^{n}\xrightarrow{}\mathbb{R}^{n}$ and $G: \mathbb{R}^{n}\rightarrow\mathbb{R}^{n}$ be smooth maps. We say that $F$ and $G$ are $k$--jet equivalent at $p \in\mathbb{R}^{n}$ if $F(p) = G(p)$ and $F(x) - G(x) = \mathcal{O}\left( ||x-p||^{k+1} \right)$ as $x\rightarrow p$. Moreover, an equivalence class defined by this concept is called the $k$--jet of $F$ at $p$, and is denoted as $j^{k}F(p)$.
        \end{definition}

       Stabilising a canard is a more delicate objective than stabilising a hyperbolic
equilibrium, and it requires a corresponding restriction on the feedback. The
fold, its contact order, and hence the entire canard family are local objects, determined
by the structure of the vector field in a neighbourhood of the fold point. A feedback that
ignores this structure could meet a naive stabilisation target by altering the system
rather than controlling it: it could, for instance, regularise the fold into a normally
hyperbolic point, or change its contact order, in either case replacing the canard one wants to stabilise by an object of a different type. To
exclude this, we require the open-- and closed--loop vector fields to share the leading
jet that defines the singularity at the fold (Definition~\ref{Def:KJetEquivalence}). This is made precise as follows.

\begin{definition}[Compatible controller]
Given a control system
\[
  \dot{\hat\zeta} = f(\hat\zeta,\hat\lambda,u),
\]
where $\hat\zeta\in\mathbb R^{n}$ denotes the state variable, $\hat\lambda\in\mathbb R^{p}$
represents the set of parameters in the system, and $u\in\mathbb R^{m}$ stands for the
control entry. Suppose that in the open--loop system, i.e.\ ($u=0$), the origin
$\hat\zeta=0$ is a nilpotent equilibrium point of $\dot{\hat\zeta}=f(\hat\zeta,0,0)$ and
that there exists a $k\in\mathbb N$ such that $k$ is the smallest number so that
$j^{k}f(0)\neq 0$. Furthermore, the control variable is $u(\hat\zeta,\hat\lambda,l)$, where
$l\in\mathbb R^{m}$ represent all the controller parameters, and let
$\dot{\hat\zeta}=F(\hat\zeta,\hat\lambda,l)$ be the closed--loop system. Then, $u$ is a
compatible controller if the open and closed--loop vector fields are $k$--jet equivalent at
the origin for $\hat\lambda=0$ \cite{Ref:JardonKojakhmetov2022}.
\end{definition}

In other words, a compatible controller is free to reshape the higher--order terms,
and thereby to select which level set $\{H=h\}$, i.e.\ which canard cycle, becomes
attracting, but it leaves the local singularity, the fold and its contact order, unchanged.
The stabilised orbit is then genuinely a canard of the original system rather than an artefact engineered by the control.

\section*{Acknoeledgements}
We sincerely thank the anonymous reviewers whose comments helped us to improve the manuscript.

\bibliographystyle{abbrv}
\bibliography{Ref}

\end{document}